\documentclass{amsart}
\usepackage{graphicx}
\usepackage{amsmath}
\usepackage{amsfonts}
\usepackage{amssymb}
\usepackage{ifpdf}
\newtheorem{theorem}{Theorem}

\newtheorem{corollary}[theorem]{Corollary}

\newtheorem{definition}[theorem]{Definition}
\newtheorem{example}[theorem]{Example}

\newtheorem{lemma}[theorem]{Lemma}

\newtheorem{proposition}[theorem]{Proposition}
\newtheorem{remark}[theorem]{Remark}

\def\Xint#1{\mathchoice
    {\XXint\displaystyle\textstyle{#1}}    {\XXint\textstyle\scriptstyle{#1}}    {\XXint\scriptstyle\scriptscriptstyle{#1}}    {\XXint\scriptscriptstyle\scriptscriptstyle{#1}}      \!\int}
\def\XXint#1#2#3{{\setbox0=\hbox{$#1{#2#3}{\int}$}
    \vcenter{\hbox{$#2#3$}}\kern-.5\wd0}}
\def\dint{\Xint-}

\makeatletter \@namedef{subjclassname@2020}{\textup{2020}
Mathematics Subject Classification} \makeatother

\begin{document}
\title[]{Generalised Haj{\l}asz-Besov spaces on $RD$-spaces }
\author{Joaquim Mart\'{i}n$^{\ast}$}
\address{Department of Mathematics\\
Universitat Aut\`{o}noma de Barcelona\\ Bellaterra, Barcelona.
Spain} \email{Joaquin.Martin@uab.cat\\\emph{ORCID:
0000-0002-7467-787X}}

\author{Walter A.  Ortiz**}
\address{School of Engineering, Science and Technology\\
Universidad Internacional de Valencia\\
Valencia. Spain\\
}
\email{walterandres.ortiz@professor.universidadviu.com\\\emph{ORCID:
0000-0002-8617-3919}}
\thanks{$^{\ast}$Partially supported by Grants PID2020-113048GB-I00, PID2020-114167GB-I00 funded  both by MCIN/AEI/10.13039/501100011033 and Grant 2021-SGR-00071
(AGAUR, Generalitat de Catalunya).}
\thanks{**Partially supported by and Grant 2021-SGR-00071
(AGAUR, Generalitat de Catalunya).}
\thanks{This paper is in final form and no version of it will be submitted for
publication elsewhere.} 

\begin{abstract}
An $RD$ space is a doubling measure metric space $\Omega$ with the
additional property that it has a reverse doubling property.  In
this paper we introduce a new class of Haj{\l}asz-Besov spaces on
$\Omega$ and extend several results from classical theory, such as
embeddings and Sobolev-type embeddings.
\end{abstract}
\maketitle

\section{INTRODUCTION}

Recently, Besov type spaces defined on metric spaces have been
rapidly introduced and developed, and several results from classical
theory have been extended into this new context (see for example
\cite{Tri0}, \cite{HMY}, \cite{MY}, \cite{Amiran}, \cite{Amiran1},
\cite{HIH}, \cite{CaoGri}, \cite{LiYaYu}, \cite{SDY}, \cite{AYY},
\cite{LS} and references therein). In this paper, we continue our
research on this topic initiated in our previous papers \cite{MaOr},
\cite{MarOrt} and \cite{MarOrt1}.

Among the various equivalent expressions for defining Besov spaces
in metric measure spaces, we use the approach based on a
generalisation of the classical $L^{p}$ modulus of smoothness
introduced in \cite{Amiran}. More specifically (precise definitions
and properties used in this introductory section can be found in
section \ref{section2}), we assume that $(\Omega,d, \mu)$ is a
metric measure space equipped with a metric $d$ and a Borel regular
outer measure $\mu$, for which the measure of each ball is positive
and finite and $\mu\left( \left\{  x\right\} \right) =0$ for all
$x\in\Omega$. Given $t>0,$ $1\leq p<\infty$ and $f\in
L_{loc}^{p}(\Omega),$ the $L^{p}$-modulus of smoothness is defined
by
\[
\mathcal{E}_{p}(f,t)=\left(  \int_{\Omega}\left(
\dint_{B(x,t)}\left| f(x)-f(y)\right|  ^{p}d\mu(y)\right)
d\mu(x)\right)  ^{1/p},
\]
where, as usual,
$\dint_{A}g(x)d\mu(y)=\frac{1}{\mu(A)}\int_{A}g(x)d\mu(y)$ is the
integral average of the function $g$ over a measurable set $A$ with
$0<\mu(A)<\infty$.

For $0<\theta<\infty,$ the homogeneous Besov space $\mathcal{\dot{B}}%
_{p,q}^{\theta}(\Omega)$ consists of functions $f\in L_{loc}^{p}(\Omega)$ for
which the seminorm
\begin{equation*}
\left\|  f\right\|  _{\mathcal{\dot{B}}_{p,q}^{\theta}(\Omega)}:=\left\{
\begin{array}
[c]{ll}%
\left(  \int_{0}^{\infty}\left(  t^{-\theta}\mathcal{E}_{p}(f,t)\right)
^{q}\frac{dt}{t}\right)  ^{1/q}, & 1\leq q<\infty,\\
\sup_{t>0}t^{-\theta}\mathcal{E}_{p}(f,t), & q=\infty,
\end{array}
\right.
\end{equation*}
is finite.

If $0<\theta<1$, this definition gives the usual Besov space in the
Euclidean setting, since as was observed in \cite{Amiran}, $\mathcal{E}%
_{p}(f,t)$ is equivalent to the standard $L^{p}$-modulus of smoothness
$\omega_{p}(f,t)=\sup_{\left|  h\right|  \leq t}\left\|  f(x+h)-f(x)\right\|
_{L^{p}(\mathbb{R}^{n})}$.

\

Classical  theory generalises Besov spaces in the following way:

\begin{enumerate}
\item  Replacing the regularity index $t^{-\theta}$ by $t^{-\theta}b(t)$,
where $b$ is a slowly varying function ($b\in SV)$.  This change
allows to modify the smoothness properties of functions; the
resulting spaces are often called generalised smooth Besov spaces,
see e.g. \cite{DoTi}, \cite{Gold01}, \cite{Gol02}, \cite{Tri00},
\cite{Tri1}, \cite{Tri2}, the literature on this subject is
extensive.

\item  Replacing the basis space $L^{p}$ in the definition of the $L^{p} $
modulus of smoothness by a Lorentz space, an Orlicz space or a
rearrangement invariant (r.i) function space\footnote{i.e. a
function space which satisfies the following condition: if $f$ and
$g$ have the same distribution function, then $\|f\|_X = \|g\|_X $
(see Section \ref{section21} below).} (see e.g. \cite{FS},
\cite{ST}, \cite{BeCo}, \cite{BCT}, \cite{Dit} and the references
given there).
\end{enumerate}

The combination of these two ideas has been a frequently used means
of obtaining the complete solution of some natural questions that
require a fine control of the behaviour of the functions involved
(for example in variational problems and PDE's), and numerous works
have been devoted to proving embedding theorems and fundamental
topological properties for Besov spaces obtained in this way, which
are crucial in applications\footnote{We refer the interested reader,
for instance, to \cite{Mami}, \cite{GPS}, \cite{Musi}, \cite{Dit},
\cite{BK}, \cite{GKe},\cite{God}, \cite{Akg}, \cite{AACS},
\cite{BOT}, \cite{SCAB}, \cite{BMO} and the references quoted
therein.}.

In this paper we will introduce a new type of Besov spaces defined
in metric measure spaces, which incorporates and generalises the
ideas described in items (i) and (ii) above. More specifically,

\begin{definition}
Let $\left(  \Omega,d,\mu\right)  $ be an $RD-$ space. Let $X$ be an
r.i. space on $\Omega,$ $E$ an r.i. space on $\left( 0,\infty\right)
$, $b\in SV $ and $0<\theta<\infty.$ The \textbf{Generalised
Homogeneous Haj\l asz-Besov Space} \r{B}$_{X,E}^{\theta,b},$ is
defined as the set of functions $f\in L^{1}+L^{\infty}$ such that
the quasi-norm
\begin{equation}
\left\|  f\right\|  _{\text{\r{B}}_{X,E}^{\theta,b}}=\left\|  t^{-\theta
}b(t)E_{X}(f,t)\right\|  _{\tilde{E}} \label{norbes}%
\end{equation}
is finite. Here $E_{X}(f,r)$ is the the $X-$modulus of continuity defined by
\[
E_{X}(f,t)=\left\|  \dint_{B(x,t)}\left|  f(x)-f(y)\right|
d\mu(y)\right\|
_{X},\text{\ }%
\]
and $\tilde{E}$ is the corresponding r.i. space with respect to
homogeneous measure $dt/t$.\footnote{See Section \ref{section221}
below.}

Similarly, the \textbf{Generalised Haj{\l}asz-Besov space} B$_{X,E}^{\theta
,b}$ is \r{B}$_{X,E}^{\theta,b}\cap X$ with the norm
\[
\left\|  f\right\|  _{\text{B}_{X,E}^{\theta,b}}=\left\|  t^{-\theta}%
b(t)E_{X}(f,t)\right\|  _{\tilde{E}(0,1)}+\left\|  f\right\|  _{X}.
\]

\end{definition}

\begin{example}
If $0<\theta<1$ and $E=L^{q}(0,\infty)$ then $\tilde{E}=L^{q}%
((0,\infty),\frac{dt}{t})$ and \r{B}$_{L^{p},L^{q}}^{\theta,1}=\mathcal{\dot{B}%
}_{p,q}^{\theta}(\Omega)$ (see Remark \ref{equal} below), while
 \r {B}$_{L^{p},L^{q}}^{\theta,b}(\mathbb{R}^{n})$ is the
classical Besov space of generalised smoothness.
\end{example}

The main point in (\ref{norbes}) is the replacement of the r.i.
space $L^{q}(0,\infty )$ by the richer family of r.i. spaces $E$,
which we believe is the natural framework for presenting this
theory. For example, any ultrasymmetric space (see \cite{Pust}) can
be derived in this way, in fact the theory of utrasymmetric spaces
developed by E. Pustylnik in \cite{Pust}, and a variant of the real
interpolation method introduced and studied in \cite{FerSig} and
\cite{FerSig1} have been both a motivation and a model for the Besov
type spaces introduced above.

The aim of this paper is to extend several results from classical
theory, such as embeddings, Sobolev-type embeddings, embedding in
BMO, essential continuity and Morrey type embeddings, for the
generalised Haj{\l}asz-Besov spaces introduced in the definition
above. Notice that Besov spaces defined in this way are also new in
the Euclidean setting; moreover, as we will see in sections
\ref{section3} and \ref{section4} below, this general context will
be of particular interest since, rather than complicating things, it
often makes the results more transparent.

The structure of this paper is as follows. In section \ref{section2}
we review well-known results and give some further notation and
background information, as well as more details on metric spaces,
function spaces, $SV $-functions and parameters. In section
\ref{section3} we prove some embedding results for generalised Haj\l
asz-Besov spaces. In section \ref{section4} we will see that for
$0<\theta<1$ generalised Haj\l asz-Besov spaces can be realised as
real interpolation spaces, this description will allow us to obtain
Sobolev-type embedding theorems. Criteria for essential continuity,
embedding in BMO and Morrey type results are also obtained and some
examples are given.

Throughout the paper we will write $f\preceq g$ instead of $f\leq
Cg$ for some constant $C>0$. The functions $f$ and $g$ are
equivalent, $f\simeq g$, if $f\preceq g$ and $g\preceq f$. We also
say that a function $f$ is almost increasing (almost decreasing) if
it is equivalent to an increasing (decreasing) function.

\section{Notation and preliminaries\label{section2}}

A metric measure space $(\Omega,d,\mu)$ is a metric space
$(\Omega,d)$ endowed with a Borel measure $\mu$ such that
$0<\mu(B)<\infty$, for every ball $B$ in $\Omega$. We will also
assume that $(\Omega,d,\mu)$ is atom free, i.e. $\mu(\{x\})=0$ for
all $x\in\Omega$. For any $x\in\Omega$ and $r>0 $, set
$B(x,r):=\{y\in\Omega:d(x,y)<r\}.$

A metric measure space $(\Omega,d,\mu)$ is called an
$RD-$\textbf{space} if there are positive constants$\
0<C_{1}\leq1\leq C_{2}$ and $0<k\leq n$ such that for all
$x\in\Omega$ and $\lambda\geq 1$, we have
\begin{equation}
C_{1}\lambda^{k}\mu(B(x,r))\leq\mu(B(x,\lambda r))\leq C_{2}\lambda^{n}%
\mu(B(x,r)). \label{iii}%
\end{equation}
We will call the pair $(k,n)$ the \textbf{indices} of the $RD$ space
(see \cite{KYZ} and \cite{YZ} for more information and equivalent
characterisations of $RD$-spaces).

We note that (\ref{iii}) implies the \textbf{doubling property},
there exists a constant $C_D\geq1$ (called the doubling constant)
such that
\begin{equation}
\mu(B(x,2r))\leq C_D\mu(B(x,r)) \label{dobla}%
\end{equation}
and the \textbf{reverse doubling property}: there exists a constant
$a\in(1,\infty)$ such that for all $x\in\Omega,\mu (B(x,ar))\geq
2\mu (B(x,r))$. \

It follows from (\ref{iii}) that there exist positive constants $c_{0},$
$C_{0}$ such that
\begin{equation*}
c_{0}\min(r^{k},r^{n})\mu(B(x,1))\leq\mu(B(x,r))\leq C_{0}\max(r^{k},r^{n}%
)\mu(B(x,1)),
\end{equation*}
for all $x\in\Omega$ and $0<r<\infty$.

We will say that a metric measure space $(\Omega,d,\mu)$ satisfies
the \textbf{non-collapsing property} if \
\begin{equation*}
\kappa=\inf_{x}\mu(B(x,1))>0.
\end{equation*}
Therefore, $RD-$spaces satisfying the non-collapsing property, verify that
\begin{equation}
c_{0}\kappa\min(r^{k},r^{n})\leq\mu(B(x,r)). \label{belw}%
\end{equation}
See \cite{MarOrt1} for more information on metric measure spaces
satisfying the non-collapsing property.

\subsection{Rearrangements of functions and rearrangement invariant spaces\label{section21}.}

In this part we will give some basic background from the theory of
rearrangements and r.i. spaces, which will be used in what follows
(for an exhaustive treatment of these topics, we refer the reader to
\cite{BS} and \cite{KPS}).

Let $\left(  \Omega,d,\mu\right)  $ be metric measure space. For
measurable functions $f:\Omega\rightarrow\mathbb{R},$ the \textbf{
decreasing rearrangement} $f_{\mu}^{*}$ of $f$ is given by
\[
f_{\mu}^{\ast}(s)=\inf\{t\geq0:\mu\{x\in{\Omega}:\left|  f(x)\right|
>t\}\leq s\},\text{ \ \ }s>0.
\]

In the following, we will omit the indices $\mu$ whenever it is
clear which measure we are working with.

A basic property of rearrangements states that
\[
\sup_{\mu(E)=t}\int_{E}\left|  f(x)\right|  d\mu=\int_{0}^{t}f^{\ast}(s)ds.
\]
Since $f^{\ast}$ is decreasing, the maximal function $f^{\ast\ast}$ of
$f^{\ast},$ defined by
\begin{equation*}
f^{\ast\ast}(t)=\frac{1}{t}\int_{0}^{t}f^{\ast}(s)ds, 
\end{equation*}
is also decreasing, and
\[
f^{\ast}\leq f^{\ast\ast}.
\]
We single out two subadditivity properties, if $f$ and $g$ are two $\mu
-$measurable functions on $\Omega$, then for $t>0$%

\begin{equation*}
\left(  f+g\right)  ^{\ast}(2t)\leq f^{\ast}(t)+g^{\ast}(t)
\end{equation*}
and%

\begin{equation*}
\left(  f+g\right)  ^{\ast\ast}(t)\leq f^{\ast\ast}(t)+g^{\ast\ast}(t).
\end{equation*}

The oscillation $O(f,\cdot)$ of $f$ is defined by
\[
O(f,t):=f^{\ast\ast}(t)-f^{\ast}(t),\text{ \ \ }0<t<\infty.
\]
Notice that
\begin{equation}
\frac{\partial}{\partial t}f^{\ast\ast}(t)=-\frac{O(f,t)}{t} \label{der2est}%
\end{equation}
and that the function $t\rightarrow tO(f,t)$ in increasing.

\

A Banach function space $X_{\mu}({\Omega})$ on $\left(  \Omega,d,\mu\right)  $
that satisfies the Fatou property is called a rearrangement invariant (r.i.)
space if, for any two measurable functions $f$ and $g$ the condition $g\in
X_{\mu}({\Omega})$, $f^{\ast}(t)\leq g^{\ast}(t)$ for almost all $t>0$ implies
$f\in X_{\mu}({\Omega})$ and $\left\|  f\right\|  _{X_{\mu}({\Omega})}%
\leq\left\|  g\right\|  _{X_{\mu}({\Omega})} $. Typical examples of r.i.
spaces are $L_{\mu}^{p}({\Omega})$-spaces, Lorentz spaces, Marcinkiewicz
spaces, Lorentz-Zygmund spaces, Orlicz spaces and Lorentz-Orlicz spaces.

To condense the notation throughout the paper, we will omit the measure and
the set when referring to a function space defined in $\left(  \Omega
,d,\mu\right)  $, so we write $X$ instead of $X_{\mu} ({\Omega})$.

For any r.i. space $X$ we have
\begin{equation*}
L^{\infty}\cap L^{1}\subset X\subset L^{1}+L^{\infty},
\end{equation*}
with continuous embedding.

A useful property of r.i. spaces states that
\begin{equation}
\int_{0}^{r}f^{\ast}(s)ds\leq\int_{0}^{r}g^{\ast}(s)ds,\text{
}\forall
r>0\Rightarrow\left\|  f\right\|  _{X}\leq\left\|  g\right\|  _{X}%
\label{Hardy}
\end{equation}
for any r.i. space $X$.

Given an r.i. space $X$, the set
\[
X^{\prime}=\left\{  f:\int_{\Omega}\left|  f(x)g(x)\right|  d\mu
(x)<\infty,\text{ \ }g\in X\right\}
\]
equipped with the norm
\begin{equation}
\left\|  f\right\|  _{X^{\prime}}=\sup_{\left\|  g\right\|  _{X}\leq1}%
\int_{\Omega}\left|  f(x)g(x)\right|  d\mu(x) \label{ass}%
\end{equation}
is called the associate space of $X$. It turns out that $X^{\prime}$
endowed with the norm given by (\ref{ass}) is an r.i. space.
Furthermore, the H\"{o}lder inequality
\begin{equation}
\int_{\Omega}\left|  f(x)g(x)\right|  d\mu(x)\leq\left\|  f\right\|
_{X}\left\|  g\right\|  _{X^{\prime}} \label{Hold}%
\end{equation}
holds for every $f\in X$ and $g\in X^{\prime}$.

A useful tool, in the study of an r.i. space $X$ is the \textbf{fundamental
function\ }of $X$ defined by
\[
\phi_{X}(t)=\left\|  \chi_{F}\right\|  _{X},
\]
where $F$ is any measurable subset of $\Omega$ with $\mu(F)=t.$ The
fundamental function $\phi_{X}$ is quasi-concave. Moreover, one has that
\begin{equation}
\phi_{X}(t)\phi_{X^{\prime}}(t)=t  \label{si}.
\end{equation}

Let $1\leq p<\infty$ and let $X$ be an r.i. space on $\Omega;$ the
$p-$\textbf{convexification} $X^{(r)}$ of $X,$ (see \cite{LT}) is the r.i
space defined \textbf{\ }by
\[
X^{(p)}=\{f:\left|  f\right|  ^{p}\in X\},\text{ \ \ }\left\|  f\right\|
_{X^{(p)}}=\left\|  \left|  f\right|  ^{p}\right\|  _{X}^{1/p}.
\]
For example $(L^1)^{(p)}=L^p$.
\subsection{Parameters spaces, extension indices and slowly varying functions\label{section22}}

\subsubsection{Parameters Spaces}

We consider two different measures on $(0,\infty),$ the usual
Lebesgue measure $dt$ and the homogeneous measure $dt/t$. We use
characters with a tilde for spaces with the measure $dt/t$. For
example,
\[
\left\|  f\right\|  _{\tilde{L}^{1}}=\int_{0}^{\infty}\left|  f(t)\right|
\frac{dt}{t}\text{ }%
\]
while $\tilde{L}^{\infty}(0,\infty)=L^{\infty}(0,\infty)$.

A \textbf{parameter} space $E=E(0,\infty)$ will be  an r.i. space on
$(0,\infty)$ with respect to the Lebesgue measure. Given a parameter
$E$, we may always assume that $E$ is an exact interpolation space
between $L^{1}$ and $L^{\infty}$ i.e.
$E:=\mathcal{F}(L^{1},L^{\infty}) $ for some real interpolation
functor $\mathcal{F}.$ Together with $E$, we consider the space
obtained by the same functor $\mathcal{F}$ from the couple
$(\tilde{L}^{1},L^{\infty})$, i.e.
$\tilde{E}:=\mathcal{F}(\tilde{L}^{1},L^{\infty}).$ The space
$\tilde{E}$ is an r.i. space respect to the measure $dt/t.$ The
spaces $E=E(0,\infty )$ and $\tilde{E}=\tilde{E}(0,\infty )$ can
also be connected directly without using interpolation functors,
since $f(t)\in
\tilde{E}(0,1)$ if and only if $f(e^{-u})\in $ $E$ and $f(t)\in \tilde{E}%
(1,\infty )$ if and only if $f(e^{u})\in E$. Moreover, $\left\| f\right\| _{%
\tilde{E}(0,1)}=\left\| f(e^{-u})\right\| _{E(0,\infty )}$ and
$\left\| f\right\| _{\tilde{E}(1,\infty )}=\left\| f(e^{u})\right\|
_{E(0,\infty )}$ (see \cite{Pust}).

Sometimes it is useful to divide the interval $(0,\infty)$ into two
subintervals $(0,t)$ and $(t,\infty)$, considering the spaces
$E(0,t)$ and $E(t,\infty)$ separately, using this notation, we can
write
\[
\left\|  f\right\|  _{E}\simeq\left\|  f\right\|  _{E(0,t)}+\left\|
f\right\|  _{E(t,\infty).}%
\]
The same applies to $\tilde{E}$.

A useful property of $\tilde{E},$ (see \cite{Pust1}), is that for
any positive numbers $a,b>0$ and any function $g\in E$, the
following equivalence holds
\begin{equation}
\left\| g(at^{b})\right\| _{\tilde{E}}\leq\min(1,b)\left\| g\right\|
_{\tilde{E}}. \label{propi}
\end{equation}

Associated with each parameter $E$ is a pair of indices, the upper
and lower and lower Boyd indices (see \cite{BS}), defined by
\begin{equation*}
\overline{\alpha}_{E}=\inf\limits_{s>1}\dfrac{\ln h_{E}(s)}{\ln
s}\text{ \ \ and \ \
}\underline{\alpha}_{E}=\sup\limits_{s<1}\dfrac{\ln h_{E}(s)}{\ln
s}.
\end{equation*}
where $h_{E}(s)$ is the norm of the dilation operator $D_{\frac{1}{s}%
}f(t)=f\left(  \frac{t}{s}\right)  ,s>0$. For example
$\underline{\alpha }_{L^{q}}=\overline{\alpha}_{L^{q}}=1/q$.

For more information on this topic see \cite{CwiPus}, \cite{Pust},
\cite{Pust1}, \cite{FerSig} and the references quoted therein.

\subsubsection{Slowly varying functions\label{section221}}

We say that a positive Lebesgue measurable function $b$, is
\textbf{slowly varying} on $(0,\infty)$ ($b\in SV)$ if, for any
$\varepsilon>0$, the function $t^{\varepsilon}b(t)$ is almost
increasing on $(0,\infty)$ while the function $t^{-\varepsilon}b(t)$
is almost decreasing on $(0,\infty)$. In fact, we have that for
every $\varepsilon>0$,  there is a positive constant
$c_{\varepsilon}$ such that
\begin{equation}
c_{\varepsilon}^{-1}\min(t^{-\varepsilon},t^{\varepsilon})b(s)\leq b(st)\leq
c_{\varepsilon}\max(t^{-\varepsilon},t^{\varepsilon})b(s),\text{ }t,s>0.
\label{sv1}%
\end{equation}

The family of $SV$ functions includes powers of logarithms, $\ell
(t)=(1+\left|  \ln t\right|  )^{\alpha}$, $t>0$ and $\alpha\in\mathbb{R},$
broken logarithmic functions
\[
\ell^{\mathbb{A}}(t)=\left\{
\begin{array}
[c]{ll}%
(1+\left|  \ln t\right|  )^{\alpha}, & 0<t<1,\\
(1+\left|  \ln t\right|  )^{\beta}, & t\geq1,
\end{array}
\right.
\]
$\mathbb{A}=(\alpha,\beta)\in\mathbb{R}^{2}$ and also functions such
as $b(t)=\exp\left(  \left|  \ln t\right|  ^{\alpha}\right)  ,$
$\alpha\in\left( 0,1\right)$. The latter family of functions has the
interesting property that each of its members tends to infinity
faster than any positive power of the logarithmic function.

The following result (see \cite{FerSig}) will be useful in what follows.

\begin{lemma}
\label{pamam}Let $E$ be a parameter and $b\in SV$. For all
$\alpha>0$ we get

\begin{enumerate}
\item
\[
\left\|  s^{\alpha}b(s)\right\|  _{\tilde{E}(0,t)}\simeq t^{\alpha}b(t)\text{
\ and }\left\|  s^{-\alpha}b(s)\right\|  _{\tilde{E}(t,\infty)}\simeq
t^{-\alpha}b(t).
\]

\item  When $\alpha>0$, the Hardy-type inequality
\[
\left\|  s^{-\alpha}b(s)\int_{0}^{t}f(s)ds\right\|  _{\tilde{E}}\leq\left\|
s^{1-\alpha}b(s)f(s)ds\right\|  _{\tilde{E}}%
\]
holds for each measurable positive function f on $(0,\infty)$.
\end{enumerate}
\end{lemma}

For more information and further examples of slowly varying functions see
\cite{BGT} or \cite{Mtz}.

\subsubsection{Extension indices}

Let $\varphi$ be a positive finite function on $(0,\infty)$, we
define its associated dilation function as
\[
M_{\varphi}(t)=\sup_{s>0}\frac{\varphi(st)}{\varphi(s)}%
\]
If $M_{\varphi}(t)$ is finite everywhere, then the\textbf{\ lower
and upper extension indices} of $\varphi$ exist and they are defined
by
\[
\underline{\beta}_{\varphi}=\lim_{t\rightarrow0}\frac{\log M_{\varphi}%
(t)}{\log t},\text{ \ }\overline{\beta}_{\varphi}=\lim_{t\rightarrow\infty
}\frac{\log M_{\varphi}(t)}{\log t}.
\]
In general,
$-\infty<\underline{\beta}_{\varphi}\leq\overline{\beta}_{\varphi
}<\infty$, but if it is increasing, then
$0\leq\underline{\beta}_{\varphi
}\leq\overline{\beta}_{\varphi}<\infty$, and if it is quasi-concave,
we have
$0\leq\underline{\beta}_{\varphi}\leq\overline{\beta}_{\varphi}\leq1$.
Also
\begin{equation}
\varphi(t)\simeq\int_{0}^{t}\varphi(s)\frac{ds}{s}\text{ if
}0<\underline {\beta}_{\varphi}\text{ and
}\varphi(t)\simeq\int_{t}^{\infty}\varphi (s)\frac{ds}{s}\text{ if
}\overline{\beta}_{\varphi}<0. \label{indi}
\end{equation}
Note also that both indices remain the same after replacing
$\varphi$ by any equivalent function.

The following properties of extension indices are easy to prove:

\begin{proposition}
\label{sv}Let $b\in SV$ and $f,g$ positive finite function on $(0,\infty)$. Then

\begin{enumerate}
\item $\overline{\beta}_{b}=\underline{\beta}_{b}=0.$

\item  If $\varphi(t)=t^{a}b(t)f(t),$ $a\in\mathbb{R}$, then $\underline
{\beta}_{\varphi}=a+\underline{\beta}_{f}$ and $\overline{\beta}_{\varphi
}=a+\overline{\beta}_{f}.$

\item  If $\varphi(t)=f(1/t)$ or $\varphi(t)=1/f(t),$ then $\underline{\beta
}_{\varphi}=-\overline{\beta}_{f}$ and $\overline{\beta}_{\varphi}%
=-\underline{\beta}_{f}.$

\item  If $\varphi(t)=f(t)g(t),$ then $\underline{\beta}_{\varphi}%
\geq\underline{\beta}_{f}$ $+\underline{\beta}_{g}$ and $\overline{\beta
}_{\varphi}\leq\overline{\beta}_{f}+\overline{\beta}_{g}.$
\end{enumerate}
\end{proposition}

It follows from (ii) that $\varphi(t)/t^{a}$ is almost increasing
for any $a<\underline{\beta}_{\varphi}$ and almost decreasing for
any $a>\overline {\beta}_{\varphi}$.

Given $X$  an r.i. space on $\Omega$, the {\bf fundamental indices}
$\underline{\beta }_{X}$ and $\bar{\beta}_{X}$ of $X$ are defined as
the lower and upper extension indices of its fundamental function
$\varphi_X$.

\section{\label{section3}Embeddings of Generalised Haj{\l}asz-Besov spaces}

In this section we first analyse the effect of replacing the basis space $X$
and/or the parameter $E$ by convexifications of them, and we prove several
embedding theorems.

\begin{proposition}{\label{emm}}
Let $\left(  \Omega,d,\mu\right)  $ be an $RD-$space, $X$ an r.i.
space, $E$ a parameter, $b\in SV$ and $0<\theta<\infty,$ then
\[
\text{\r{B}}_{X,L^{1}}^{\theta,b}\subset\text{\r{B}}_{X,E}^{\theta,b}%
\subset\text{\r{B}}_{X,L^{\infty}}^{\theta,b}%
\]
and
\[
B_{X,L^{1}}^{\theta,b}\subset B_{X,E}^{\theta,b}\subset B_{X,L^{\infty}%
}^{\theta,b}.
\]
\end{proposition}

\begin{proof}
It follows from (\ref{dobla}) that there is a positive constant $c$ which
depends only on the doubling constant, such that
\begin{equation}
c^{-1}E_{X}(f,s)\leq E_{X}(f,R)\leq cE_{X}(f,r) \label{qc}
\end{equation} for
all $0<R/2\leq s\leq R\leq r\leq2R$, thus
\begin{align}
b(r)r^{-\theta}E_{X}(f,r)  &  \preceq\int_{r}^{2r}b(t)t^{-\theta}%
E_{X}(f,t)\frac{dt}{t}\leq\left\|  b(t)t^{-\theta}E_{X}(f,t)\right\|
_{\tilde{E}}\left\|  \chi_{(r,2r)}\right\|  _{\tilde{E}^{\prime}}\label{cal}\\
&  =\left\|  f\right\|  _{\text{\r{B}}_{X,E}^{\theta,b}}\varphi_{\tilde
{E}^{\prime}}(\ln2).\nonumber
\end{align}
Therefore
\[
\text{\r{B}}_{X,E}^{\theta,b}\subset\text{\r{B}}_{X,L^{\infty}}^{\theta,b}.
\]

To see the first embedding, we will see that for any $t>0$ we have
\[
\left\|  b(r)r^{-\theta}E_{X}(f,r)\right\|  _{\tilde{E}(0,t)}\preceq\int
_{0}^{2t}b(s)s^{-\theta}E_{X}(f,s)\frac{ds}{s}.
\]
In virtue of interpolation properties of any parameter $E$, it
suffices to prove this inequality for the extreme spaces
$\tilde{L}^{1}$ and $L^{\infty}$. The above inequality is obvious
for $\tilde{E}=\tilde{L}^{1}$. Let us check it for $\tilde{E}$
$=L^{\infty}$.

By (\ref{cal}) we get
\[
b(r)r^{-\theta}E_{X}(f,r)\preceq\int_{0}^{2r}b(s)t^{-\theta}E_{X}%
(f,s)\frac{ds}{s},
\]
thus
\[
\left\|  b(r)r^{-\theta}E_{X}(f,r)\right\|  _{L^{\infty}(0,t)}\preceq\int
_{0}^{2t}b(s)s^{-\theta}E_{X}(f,s)\frac{ds}{s}.
\]
\end{proof}

The next result is an extension of the well-known embedding $B_{p,q_{1}%
}^{\theta}\subset B_{p,q_{2}}^{\theta}$ when $1\leq q_{1}<q_{2}$.

\begin{proposition}
Let $\left(  \Omega,d,\mu\right)$ be an $RD-$space, $X$ an r.i.
space on $\Omega$, $E$ a parameter, $b\in SV$ and $0<\theta<\infty$.
Suppose that $1\leq p<q<\infty$, then
\[
\left\|  f\right\|
_{\text{\r{B}}_{X,E^{(q)}}^{\theta,b}}\preceq\left\| f\right\|
_{\text{\r{B}}_{X,E^{(p)}}^{\theta,b}}.
\]
The same embedding holds for the inhomogeneous case.
\end{proposition}

\begin{proof}
By the previous result, \r{B}$_{X,E^{(p)}}^{\theta,b}\subset$\r{B}%
$_{X,L^{\infty}}^{\theta,b},$ hence
\begin{align*}
\left\|  f\right\|  _{\text{\r{B}}_{X,E^{(q)}}^{\theta,b}}^{q}  &  =\left\|
\left(  t^{-\theta}b(t)E_{X}(f,t)\right)  ^{q}\right\|  _{\tilde{E}}=\left\|
\left(  t^{-\theta}b(t)E_{X}(f,t)\right)  ^{q-p}\left(  t^{-\theta}%
b(t)E_{X}(f,t)\right)  ^{p}\right\|  _{\tilde{E}}\\
&  \leq\sup_{t}\left(  t^{-\theta}b(t)E_{X}(f,t)\right)  ^{q-p}\left\|
\left(  t^{-\theta}b(t)E_{X}(f,t)\right)  ^{p}\right\|  _{\tilde{E}}\\
&  =\left\|  f\right\|  _{B_{X,L^{\infty}}^{\theta,b}}^{q-p}\left\|
t^{-\theta}b(t)E_{X}(f,t)\right\|  _{\tilde{E}^{(p)}}^{p}\\
&  \preceq\left\|  t^{-\theta}b(t)E_{X}(f,t)\right\|  _{\tilde{E}^{(p)}}^{q},
\end{align*}
that implies the desired result.
\end{proof}

We will now consider the effect of replacing the space $X$ by
$X^{(q)}$ and the relation between the spaces \r{B}$_{X,E}^{\theta
,b}$ and \r{B}$_{X^{(q)},E}^{\theta ,b}$, for which we need the
following technical result obtained by Ranjbar-Motlagh in \cite{Ran}
for the special case $X=L^{q}$. We will give its proof in the
Appendix section at the end of this paper.

\begin{lemma}
\label{Suerte}Let $\left(  \Omega,d,\mu\right)  $ be an $RD-$space
that satisfies the non-collapsing condition, $X$ an r.i. space on
$\Omega$, $q\geq p\geq1$. Then, for any $f\in L^{1}+L^{\infty}$ and
any $r>0$, there exists a constant $c=c(c_{0},\kappa,C_{D})$ such
that

\begin{enumerate}
\item
\[
\left\|  \dint_{B(x,r)}\left|  f(y)\right|  d\mu(y)\right\|  _{X^{(q)}}%
\leq\frac{c}{\varphi_{X}\left(  \min(r^{k},r^{n})\right)  ^{\frac{1}{p}%
-\frac{1}{q}}}\left\|  f\right\|  _{X^{(p)}}%
\]

\item
\[
E_{X^{(q)}}(f,r)\leq c\int_{0}^{4r}\frac{\mathcal{E}_{X^{(p)}}(f,s)}%
{\varphi_{X}\left(  \min(s^{k},s^{n})\right)  ^{\frac{1}{p}-\frac{1}{q}}}%
\frac{ds}{s},
\]
where
\[
\mathcal{E}_{X^{(p)}}(f,t):=\left\|  \left(  \dint_{B(x,t)}\left|
f(x)-f(y)\right|  ^{p}d\mu(y)\right)  ^{1/p}\right\|  _{X^{(p)}}.
\]
\end{enumerate}
\end{lemma}

\begin{theorem}
{\label{anterior}} Let $\left( \Omega ,d,\mu \right) $ be an
$RD-$space, $X$ an r.i. space, $E$ a parameter, $b\in SV,$ $0<\theta
<\infty ,$ $1\leq p<\infty \ $and $1<\sigma <\infty .$ Let
$\bar{\beta}_{X}\leq \gamma \leq 1$ be such that the function
$\frac{\phi _{X}(t)}{t^{\gamma }}$ is quasi
decreasing\footnote{By proposition \ref{sv}, this is always satisfied if $\bar{\beta}_{X}<\gamma \leq 1$.} and $\alpha =\theta +\frac{n\gamma }{p}\left( 1-\frac{1}{\sigma }%
\right) .$ Then

\begin{enumerate}
\item
\begin{equation*}
\left\| f\right\| _{\text{\r{B}}_{X^{(\sigma
p)},E}^{\theta,b}}\preceq\left\|
f\right\|_{B_{X^{(p)},E}^{\alpha,b}},
\end{equation*}

\item
\begin{equation*}
\left\| f\right\| _{B_{X^{(\sigma p)},E}^{\theta ,b}}\preceq \left\|
f\right\| _{B_{X^{(p)},E}^{\alpha ,b}},
\end{equation*}

\item  and
\begin{equation}
\left\| f\right\| _{X^{(p\sigma )}}\preceq \left\| f\right\| _{X^{(p)}}^{%
\frac{\theta }{\alpha }}\left\| f\right\| _{B_{X^{(p)},E}^{\alpha ,b}}^{1-%
\frac{\theta }{\alpha }}\left( 1+\frac{1}{b\left( \left( \left\|
f\right\|
_{X^{(p)}}/\left\| f\right\| _{B_{X^{(p)},E}^{\alpha ,b}}\right) ^{\frac{%
\theta }{\alpha }}\right) }\right) .  \label{submult}
\end{equation}
\end{enumerate}
\end{theorem}

\begin{proof}
(i) By the previous Lemma, and since $\phi _{X}(t)/t^{\gamma }$ is
quasi decreasing, we get
\begin{align*}
\left\| f\right\| _{\text{\r{B}}_{X^{(\sigma p)},E}^{\theta ,b}}&
=\left\|
b(t)t^{-\theta }E_{X^{(\sigma p)}}(f,t)\right\| _{\tilde{E}} \\
& \preceq \left\| b(t)t^{-\theta }\int_{0}^{4t}\frac{\mathcal{E}%
_{X^{(p)}}(f,s)}{\varphi _{X}\left( \min (s^{k},s^{n})\right) ^{\frac{1}{p}-%
\frac{1}{\sigma p}}}\frac{ds}{s}\right\| _{\tilde{E}} \\
& \simeq \left\| b(4t)\left( 4t\right) ^{-\theta }\int_{0}^{4t}\frac{%
\mathcal{E}_{X^{(p)}}(f,s)}{\varphi _{X}\left( \min (s^{k},s^{n})\right) ^{%
\frac{1}{p}-\frac{1}{\sigma p}}}\frac{ds}{s}\right\| _{\tilde{E}} \\
& \simeq \left\| b(t)t^{-\theta }\int_{0}^{t}\frac{\mathcal{E}_{X^{(p)}}(f,s)%
}{\varphi _{X}\left( \min (s^{k},s^{n})\right)
^{\frac{1}{p}-\frac{1}{\sigma
p}}}\frac{ds}{s}\right\| _{\tilde{E}}\text{ (by (\ref{propi}))} \\
& \preceq \left\| b(t)t^{-\theta
}\frac{\mathcal{E}_{X^{(p)}}(f,t)}{\varphi
_{X}\left( \min (t^{k},t^{n})\right) ^{\frac{1}{p}-\frac{1}{\sigma p}}}%
\right\| _{\tilde{E}}\text{ (by Lemma \ref{pamam})} \\
& \preceq \left\| b(t)t^{-\theta }\left( \frac{t^{n\gamma }}{\varphi
_{X}\left( t^{n}\right) }\right) ^{\frac{1}{p}-\frac{1}{\sigma p}}\frac{%
\mathcal{E}_{X^{(p)}}(f,t)}{t^{n\gamma \left( \frac{1}{p}-\frac{1}{\sigma p}%
\right) }}\right\| _{\tilde{E}(0,1)}+\left\| b(t)t^{-\theta }\frac{\mathcal{E%
}_{X^{(p)}}(f,t)}{\varphi _{X}\left( t^{k}\right) ^{\frac{1}{p}-\frac{1}{%
\sigma p}}}\right\| _{\tilde{E}(1,\infty )} \\
& \preceq \left\| b(t)t^{-\theta }\frac{\mathcal{E}_{X^{(p)}}(f,t)}{%
t^{n\gamma \left( \frac{1}{p}-\frac{1}{\sigma p}\right) }}\right\| _{\tilde{E%
}(0,1)}+\left\| b(t)t^{-\theta }\mathcal{E}_{X^{(p)}}(f,t)\right\| _{\tilde{E%
}(1,\infty )}\text{  } \\
& \preceq \left\| f\right\| _{\text{\r{B}}_{X^{(p)},E}^{\alpha
,b}}+\left\| b(t)t^{-\theta }\mathcal{E}_{X^{(p)}}(f,t)\right\|
_{\tilde{E}(1,\infty )}.
\end{align*}
Finally, since
\begin{eqnarray*}
\mathcal{E}_{X^{(p)}}(f,t) &=&\left\| \left( \dint_{B(x,t)}\left|
f(x)-f(y)\right| ^{p}d\mu (y)\right) ^{1/p}\right\| _{X^{(p)}} \\
&\leq &\left\| f\right\| _{X^{(p)}}+\left\| \left(
\dint_{B(x,t)}\left|
f(y)\right| ^{p}d\mu (y)\right) ^{1/p}\right\| _{X^{(p)}} \\
&\leq &\left\| f\right\| _{X^{(p)}}+\frac{c}{\varphi _{X}(t^{k},t^{n})}%
\left\| f\right\| _{X^{(p)}\text{ }}\text{(by Lemma \ref{Suerte})}
\end{eqnarray*}
it follows that
\begin{eqnarray}
\left\| b(t)t^{-\theta }\mathcal{E}_{X^{(p)}}(f,t)\right\| _{\tilde{E}%
(1,\infty )} &\leq &\left\| b(t)t^{-\theta }\left(
1+\frac{c}{\varphi
_{X}(t^{k},t^{n})}\right) \left\| f\right\| _{X^{(p)}}\right\| _{\tilde{E}%
(1,\infty )}  \label{ultim} \\
&\leq &\left( 1+c\right) \left\| f\right\| _{X^{(p)}}\left\|
b(t)t^{-\theta
}\right\| _{\tilde{E}(1,\infty )}  \notag \\
&\preceq &\left\| f\right\| _{X^{(p)}}\text{ (by Lemma
\ref{pamam})}. \notag
\end{eqnarray}
(ii)We know from part (i)
\begin{equation*}
\left\| f\right\| _{\text{\r{B}}_{X^{(\sigma p)},E}^{\theta
,b}}\preceq \left\| f\right\| _{B_{X^{(p)},E}^{\theta ,b}}.
\end{equation*}
On the other hand
\begin{align}
\left\| f\right\| _{X^{(\sigma p)}}& \leq \left| \left\| f\right\|
_{X^{(\sigma p)}}-\left\| \dint_{B(x,1/4)}\left| f(x)\right| d\mu
(x)\right\|
_{X^{(\sigma p)}}\right|   \label{cot} \\
& +\left\| \dint_{B(x,1/4)}\left| f(x)\right| d\mu (x)\right\|
_{X^{(\sigma
p)}}  \notag \\
& =I+II.  \notag
\end{align}
From Lemma \ref{Suerte} (i) with $r=1/4$ we get
\begin{equation}
II\preceq \left\| f\right\| _{X^{(p)}}.  \label{cotII}
\end{equation}
By H\"{o}lder's inequality, the concavity of $\varphi _{X}$ and the
fact that $\phi _{X}(t)/t^{\gamma }$ is quasi decreasing, we have
\begin{align}
I& \leq \left\| \left| f(x)-\dint_{B(x,1/4)}f(y)d\mu (y)\right|
^{\sigma
p}\right\| _{X}^{1/\sigma p}  \label{cotI} \\
& \leq \left\| \left( \dint_{B(x,1/4)}\left| f(x)-f(y)\right| d\mu
(y)\right)
^{\sigma p}\right\| _{X}^{1/\sigma p}  \notag \\
& =E_{X^{(\sigma p)}}(f,1/4)\preceq \int_{0}^{1}\frac{\mathcal{E}%
_{X^{(p)}}(f,s)}{\left( \varphi _{X}\left( s^{n}\right) \right) ^{\frac{1}{p}%
-\frac{1}{\sigma p}}}\frac{ds}{s}\text{ (by Lemma \ref{Suerte}
(ii))}  \notag
\\
& =\int_{0}^{1}\frac{s^{-\theta
}b(s)\mathcal{E}_{X^{(p)}}(f,s)}{\left(
\varphi _{X}\left( s^{n}\right) \right) ^{\frac{1}{p}-\frac{1}{\sigma p}}}%
\frac{s^{\theta }}{b(s)}\frac{ds}{s}  \notag \\
& \preceq \left\| b(t)t^{-\theta }\frac{\mathcal{E}_{X^{(p)}}(f,t)}{%
t^{n\gamma \left( \frac{1}{p}-\frac{1}{\sigma p}\right) }}\right\| _{\tilde{E%
}(0,1)}\left\| \frac{s^{\theta }s^{n\gamma \left( \frac{1}{p}-\frac{1}{%
\sigma p}\right) }}{b(s)\left( \varphi _{X}\left( s^{n}\right) \right) ^{%
\frac{1}{p}-\frac{1}{\sigma p}}}\right\| _{\tilde{E}^{\prime }(0,1)}
\notag
\\
& \preceq \left\| b(t)t^{-\theta }\frac{\mathcal{E}_{X^{(p)}}(f,t)}{%
t^{n\gamma \left( \frac{1}{p}-\frac{1}{\sigma p}\right) }}\right\| _{\tilde{E%
}(0,1)}\left\| \frac{s^{\theta }}{b(s)}\right\| _{\tilde{E}^{\prime
}(0,1)}
\notag \\
& \preceq \left\| b(t)t^{-\theta }\frac{\mathcal{E}_{X^{(p)}}(f,t)}{%
t^{n\gamma \left( \frac{1}{p}-\frac{1}{\sigma p}\right) }}\right\| _{\tilde{E%
}(0,1)}\text{ (by Lemma \ref{pamam}).}  \notag
\end{align}
Inserting (\ref{cotII}) and (\ref{cotI}) into (\ref{cot}), we find
that
\begin{equation*}
\left\| f\right\| _{X^{(\sigma p)}}\leq \left\| b(t)t^{-\theta }\frac{%
\mathcal{E}_{X^{(p)}}(f,t)}{t^{\gamma n\left( \frac{1}{p}-\frac{1}{\sigma p}%
\right) }}\right\| _{\tilde{E}(0,1)}+\left\| f\right\|
_{X^{(p)}}=\left\| f\right\| _{_{B_{X^{(p)},E}^{\alpha ,b}}}.
\end{equation*}
Thus
\begin{align*}
\left\| f\right\| _{B_{X^{(\sigma p)},E}^{\theta ,b}}& =\left\|
f\right\| _{X^{(\sigma p)}}+\left\| b(t)t^{-\theta }E_{X^{(\sigma p)}}(f,t)\right\| _{%
\tilde{E}(0,1)} \\
& \preceq \left\| f\right\| _{_{B_{X^{(p)},E}^{\alpha ,b}}}
\end{align*}
as we wanted to show.

(iii) Let $0<R<1/4$, then
\begin{align*}
\left\| f\right\| _{X^{(p\sigma )}}& \leq \left| \left\| f\right\|
_{X^{(p\sigma )}}-\left\| \dint_{B(x,R)}\left| f(x)\right| d\mu
(x)\right\|
_{X^{(p\sigma )}}\right|  \\
& +\left\| \dint_{B(x,R)}\left| f(x)\right| d\mu (x)\right\|
_{X^{(p\sigma )}}
\\
& =I+II.
\end{align*}
From Lemma \ref{Suerte}$,$ we get
\begin{equation*}
II\preceq \frac{1}{\left( \varphi _{X}(R^{n})\right) ^{\frac{1}{p}-\frac{1}{%
p\sigma }}}\left\| f\right\| _{X^{(p)}},
\end{equation*}
and with the same argument that was used to get (\ref{cotI}), we
have
\begin{equation*}
I\preceq \int_{0}^{4R}\frac{\mathcal{E}_{X^{(p)}}(f,s)}{\left(
\varphi _{X}(s^{n})\right) ^{\frac{1}{p}-\frac{1}{p\sigma
}}}\frac{ds}{s}.
\end{equation*}
On the other hand, using the fact that $\varphi _{X}(t)/t^{\gamma }$
is quasi decreasing and Lemma \ref{pamam}, we have
\begin{align*}
\int_{0}^{4R}\frac{\mathcal{E}_{X^{(p)}}(f,s)}{\left( \varphi
_{X}(s^{n})\right) ^{\frac{1}{p}-\frac{1}{p\sigma }}}\frac{ds}{s}&
=\int_{0}^{4R}\frac{\mathcal{E}_{X^{(p)}}(f,s)}{s^{n\gamma }{}^{\left( \frac{%
1}{p}-\frac{1}{p\sigma }\right) }}\frac{s^{n\gamma }{}^{\left( \frac{1}{p}-%
\frac{1}{p\sigma }\right) }}{\left( \varphi _{X}(s^{n})\right) ^{\frac{1}{p}-%
\frac{1}{p\sigma }}}\frac{ds}{s} \\
& \preceq \left( \frac{R^{\gamma n}{}}{\varphi _{X}(R^{n})}\right) ^{\frac{1%
}{p}-\frac{1}{p\sigma }}\int_{0}^{4R}\frac{\mathcal{E}_{X^{(p)}}(f,s)}{%
s^{n\gamma }{}^{\left( \frac{1}{p}-\frac{1}{p\sigma }\right)
}}\frac{ds}{s}
\\
& \leq \left( \frac{R^{\gamma n}{}}{\varphi _{X}(R^{n})}\right) ^{\frac{1}{p}%
-\frac{1}{p\sigma }}\left\| \frac{b(s)s^{-\theta }\mathcal{E}_{X^{(p)}}(f,s)%
}{\min (s^{k},s^{n})^{\gamma \left( \frac{1}{p}-\frac{1}{p\sigma }\right) }}%
\right\| _{\tilde{E}(0,4R)}\left\| \frac{s^{\theta }}{b(s)}\right\| _{\tilde{%
E}^{\prime }(0,4R)} \\
& \preceq \left( \frac{R^{\gamma n}}{\varphi _{X}(R^{n})}\right) ^{\frac{1}{p%
}-\frac{1}{p\sigma }}\frac{R^{\theta }}{b(R)}\left\| \frac{b(s)s^{-\theta }%
\mathcal{E}_{X^{(p)}}(f,s)}{\min (s^{k},s^{n})^{\gamma \left( \frac{1}{p}-%
\frac{1}{p\sigma }\right) }}\right\| _{\tilde{E}}.
\end{align*}
In summary
\begin{equation*}
\left\| f\right\| _{X^{(p\sigma )}}\preceq \left( \frac{R^{\gamma n}{}}{%
\varphi _{X}(R^{n})}\right) ^{\frac{1}{p}-\frac{1}{p\sigma }}\left( \frac{1}{%
R^{n\gamma \left( \frac{1}{p}-\frac{1}{p\sigma }\right) }}\left\|
f\right\| _{X^{(p)}}+\frac{R^{\theta }}{b(R)}\left\| f\right\|
_{B_{X^{(p)},E}^{\alpha ,b}}\right) .
\end{equation*}
So for $0<R<1/4,$ one gets
\begin{equation*}
\left\| f\right\| _{X^{(p\sigma )}}\preceq \inf_{0<R<1}\left( \frac{1}{%
R^{n\gamma \left( \frac{1}{p}-\frac{1}{p\sigma }\right) }}\left\|
f\right\| _{X^{(p)}}+\frac{R^{\theta }}{b(R)}\left\| f\right\|
_{B_{X^{(p)},E}^{\alpha ,b}}\right) .
\end{equation*}
If we choose
\begin{equation*}
R=\left( \frac{\left\| f\right\| _{X^{(p)}}}{4\left\| f\right\|
_{B_{X^{(p)},E}^{\alpha ,b}}}\right) ^{\frac{1}{\alpha }},
\end{equation*}
then the submultiplicative inequality (\ref{submult}) holds.
\end{proof}

\begin{remark}
For the special case $X=L^{1}\left( \Omega\right) \ $and
$E=L^{q}(0,\infty)$, the previous theorem was proved in
\cite[Theorem 3.3]{Ran1} and \cite[Theorem 3.6]{Ran1}. See
\cite[Theorem 7.1]{KN} in the Euclidean case.
\end{remark}

\begin{example}
If $r>p\geq1$ and $\alpha=\theta+\frac{n}{p}\left(
1-\frac{p}{r}\right)  ,$ then the Theorem \ref{anterior} gives the
following submultiplicative inequality for classical Besov spaces of
logarithmic smoothness $B_{p,q}^{\theta,\ell}(\mathbb{R}^{n})$:
\[
\left\|  f\right\|  _{L^{r}(\mathbb{R}^{n})}\preceq\left\|  f\right\|
_{L^{p}(\mathbb{R}^{n})}^{\frac{\theta}{\alpha}}\left\|  f\right\|
_{B_{p,q}^{\theta,\ell}(\mathbb{R}^{n})}^{1-\frac{\theta}{\alpha}}\left(
1+\frac{1}{\left(  1+\ln\left(  \frac{\left\|  f\right\|  _{B_{p,q}%
^{\theta,\ell}(\mathbb{R}^{n})}}{\left\|  f\right\|  _{L^{p(\mathbb{R}^{n})}}%
}\right)  \right)  ^{\frac{\theta}{\alpha}}}\right).
\]
\end{example}

\section{\label{section4}Generalized Haj{\l}asz-Besov spaces as interpolation spaces}

In this section we prove that the spaces $\text{\r{B}}_{X^{(\sigma
p)},E}^{\theta,b}$ and $B_{X,E}^{\theta,b}$, $0<\theta<1$, can be
constructed using the real interpolation method described in
\cite{FerSig}; this description will allow us to obtain Sobolev-type
embedding theorems in the subsection \ref{subsection41}.

Let $\left( \Omega,d,\mu\right)$ be an $RD-$space, let $f$ be a
$\mu-$measurable function, a non-negative measurable function $g$
that satisfies
\begin{equation*}
\left|  f(x)-f(y)\right|  \leq d(x,y)\left(  g(x)+g(y)\right)  \text{ \ \ }%
\mu-a.e.\text{ \ }x,y\in\Omega.
\end{equation*}
is called a $1-$gradient of $f$. We denote by $D(f)$ the collection
of all $1-$gradients of $f$.

The \textbf{homogeneous Haj\l asz-Sobolev} space \r{M}$^{1,X}$
consists of measurable functions $f\in L^{1}+L^{\infty}$ for which
\[
\left\|  f\right\|  _{\text{\r{M}}^{1,X}}=\inf_{g\in D(f)}\left\|  g\right\|
_{X}%
\]
is finite.

The \textbf{Haj\l asz-Sobolev} space $M^{1,X}$ is \r{M}$^{1,X}\cap X$ equipped
with the quasi-norm
\[
\left\|  f\right\|  _{M^{1,X}}=\left\|  f\right\|  _{X}+\left\|  f\right\|
_{\text{\r{M}}^{1,X}}.
\]

\begin{remark}
The definition formulated above is motivated by the Haj\l asz approach to the
definition of Sobolev spaces on a metric measure space (see \cite{Ha2} and
\cite{Ha1}), where $M^{1,p}$ was defined as the set of measurable functions
$f$ for which
\[
\left\|  f\right\|  _{\text{\r{M}}^{1,p}}=\inf_{g\in D(f)}\left\|
g\right\| _{L^{p}}
\]
is finite\footnote{For $p>1$, see \cite{Ha2}, \cite{Ha1},
$M^{1,p}(\mathbb{R}^{n})=W^{1,p}(\mathbb{R} ^{n})$, while for
$n/(n+1)<p\leq1$, $M^{1,p} (\mathbb{R}^{n})$ coincides with the
Hardy-Sobolev space $H^{1,p} (\mathbb{R}^{n})$ \cite[Theorem
1]{KoSak}.}. Based on this definition \r {M}$^{1,X}$ appears
naturally when the Lebesgue norm is replaced by the norm $\left\|
\cdot\right\|  _{X}$, for example, if $X$ is an Orlicz space, we get
the Orlicz-Haj\l asz-Sobolev space (see \cite{Tuo}), and if $X$ is a
Lorentz space, we get the Lorentz-Haj\l asz-Sobolev space (see
\cite{CosMi}).
\end{remark}

Recall that given $\left(  X_{0},X_{1}\right)  $ a compatible couple of
(quasi-semi)normed spaces (i.e. there is a topological vector space
$\mathcal{X} $ such that $X_{0}$ and $X_{1}$ are continuously embedded into
$\mathcal{X)}$. For every $f\in X_{0}+X_{1}$ and $t>0$, the Peetre
$K-$functional is
\[
K(f,t;X_{0},X_{1}):=\inf_{f=f_{0}+f_{1}}\left\{  \left\|  f_{0}\right\|
_{X_{0}}+t\left\|  f_{1}\right\|  _{X_{1}}\right\}  .
\]

The following interpolation method has been introduced in
\cite{FerSig} and is an extension of the well-known real
interpolation method with a functional parameter (see \cite{Gu} and
\cite{Me}).

\begin{definition}
Let $(X_{0},X_{1})$ be a compatible couple of (quasi-semi)normed spaces, $E$ a
parameter, $b\in SV$ and $0<\theta<1$. The real interpolation space
$(X_{0},X_{1})_{\theta,b,E}$ consists of all f in $X_{0}+X_{1}$ for which the
norm
\[
\left\|  f\right\|  _{(X_{0},X_{1})_{\theta,b,E}}=\left\|  t^{-\theta
}b(t)K(f,t;X_{0},X_{1})\right\|  _{\tilde{E}}%
\]
is finite.
\end{definition}

If $E=L^{q}$ and $b=1$, then the space $(X_{0},X_{1})_{\theta,b,E}$
coincides with the classical real interpolation space
$(X_{0},X_{1})_{\theta,q}$. For more information about the spaces
$(X_{0},X_{1})_{\theta,b,E}$ and their applications, see for example
\cite{FerSig}, \cite{Me} and \cite{Per}.

\

\begin{theorem}
\label{TeoInterpol}Let $\left( \Omega ,d,\mu \right) $ be an $RD$
space, $X$
be an r.i. space, and $p\geq 1$, then for all $f\in X^{(p)}+$\r{M}$%
^{1,X^{(p)}}$ and $t>0$, we have that
\begin{align}
E_{X^{(p)}}(f,t)& \leq \mathcal{E}_{X^{(p)}}(f,t)\preceq K(f,t,X^{(p)},\text{%
\r{M}}^{1,X^{(p)}})\preceq \sum_{j=0}^{\infty
}2^{-j}E_{X^{(p)}}(f,2^{j}t)
\label{kf1} \\
& \leq \sum_{j=0}^{\infty }2^{-j}\mathcal{E}_{X^{(p)}}(f,2^{j}t)
\notag
\end{align}
and
\begin{equation*}
K(f,t,X^{(p)},M^{1,X^{(p)}})\simeq K(f,t,X^{(p)},\text{\r{M}}%
^{1,X^{(p)}})+\min (1,t)\left\| f\right\| _{X^{(p)}}.
\end{equation*}
Where
\begin{equation*}
\mathcal{E}_{X^{(p)}}(f,t):=\left\| \left(\dint_{B(x,t)}\left|
f(x)-f(y)\right| ^{p}d\mu (y)\right) ^{1/p}\right\| _{X^{(p)}}.
\end{equation*}
\end{theorem}

\begin{proof}
In \cite[Theorem 11]{MarOrt1} we proved that
\begin{equation*}
E_{X^{(p)}}(f,t)\preceq
K(f,t,X^{(p)},\text{\r{M}}^{1,X^{(p)}})\preceq \sum_{j=0}^{\infty
}2^{-j}E_{X^{(p)}}(f,2^{j}t).
\end{equation*}
By H\"{o}lder's inequality we get $E_{X^{(p)}}(f,t)\leq \mathcal{E}%
_{X^{(p)}}(f,t)$, so the right hand side inequality of (\ref{kf1})
follows. To see the left hand side inequality, given $f\in
L^{1}+L^{\infty }$, $1\leq p<\infty $ and $t>0$ we define define
\begin{equation*}
T_{r}^{p}f(x)=\left(\dint_{B(x,r)}\left| f(y)\right| ^{p}d\mu
(y)\right) ^{1/p}
\end{equation*}
we claim that the family of operators $\left\{ T_{r}^{p}\right\}
_{r>0}$ is uniformly bounded on $X^{(p)}$.

Obviously
\begin{equation}
\left\| T_{r}^{p}f\right\| _{L^{\infty \text{ }}}\leq \left\|
f\right\| _{L^{\infty \text{ }}}.  \label{a1}
\end{equation}
On the other hand, by Fubini's theorem, we get
\begin{eqnarray*}
\left\| T_{r}^{p}f\right\| _{L^{p}}^{p} &=&\int_{\Omega }\left(
\dint_{B(x,r)}\left|
f(y)\right| ^{p}d\mu (y)\right) d\mu (x) \\
&=&\int_{\Omega }\left| f(y)\right| ^{p}\left(
\int_{B(y,r)}\frac{1}{\mu (B(x,r))}d\mu (x)\right) d\mu (y).
\end{eqnarray*}
Using the doubling property and the fact that $B(x,r)\subset
B(y,2r)$ whenever $y\in B(x,r)$, we conclude that
\begin{eqnarray*}
\int_{B(y,r)}\frac{1}{\mu (B(x,r))}d\mu (x) &\leq &C_{\mu }\int_{B(y,r)}%
\frac{1}{\mu (B(x,2r))}d\mu (x) \\
&\leq &C_{\mu }\int_{B(y,r)}\frac{1}{\mu (B(y,r))}d\mu (x) \\
&=&C_{\mu }.
\end{eqnarray*}
Thus, for all $r>0$
\begin{equation}
\left\| T_{r}^{p}f\right\| _{L^{p}}^{p}\leq C_{\mu }\int_{\Omega
}\left| f(y)\right| ^{p}d\mu (y)=C_{\mu }\left\| f\right\|
_{L^{p}}^{p}.  \label{a2}
\end{equation}
By combining (\ref{a1}) and (\ref{a2}) with the definition of the
$K$ functional, we obtain
\begin{equation*}
K(T_{r}^{p}f,t^{1/p},L^{p},L^{\infty })\preceq
K(f,t^{1/p},L^{p},L^{\infty }),\text{ \ \ }(t>0).
\end{equation*}
Since (see \cite[Theorem 5.1]{BL})
\begin{equation*}
K(g,t^{1/p},L^{p},L^{\infty })\simeq \left( \int_{0}^{t}\left(
g^{\ast }(s)\right) ^{p}ds\right) ^{1/p}
\end{equation*}
we have
\begin{equation*}
\int_{0}^{t}\left( \left( T_{r}^{p}f\right) ^{\ast }(s)\right)
^{p}ds\preceq \int_{0}^{t}\left( f^{\ast }(s)\right) ^{p}ds
\end{equation*}
which by (\ref{Hardy}) implies that
\begin{equation*}
\left\| (T_{r}^{p}f)^{p}\right\| _{X}\preceq \left\| \left| f\right|
^{p}\right\| _{X},  
\end{equation*}
which is equivalent to
\begin{equation}
\left\| T_{r}^{p}f\right\| _{X^{(p)}}\preceq \left\| f\right\|
_{X^{(p)}}.  \label{asi}
\end{equation}

Let $f=g+h$, where $g\in X^{(p)}$, $h\in \dot{M}^{1,X^{(p)}}$ and
let $t>0$. \ Then
\begin{equation*}
\left(\dint_{B(x,t)}\left| g(x)-g(y)\right| ^{p}d\mu (y)\right)
^{1/p}\preceq \left| g(x)\right| +T_{r}^{p}g(x),
\end{equation*}
consequently, by the previous claim,
\begin{equation*}
E_{X^{(p)}}(g,t)\preceq \left\| g\right\| _{X^{(p)}}+\left\|
T_{r}^{p}g\right\| _{X^{(p)}}\preceq \left\| g\right\| _{X^{(p)}}.
\end{equation*}
On the other hand, since $h\in \dot{M}^{1,X^{(p)}},$ by the
definition of the $1-$gradient, if $\varrho \in D(h)\cap X^{(p)},$
then
\begin{align*}
\left(\dint_{B(x,t)}|h(x)-h(y)|^{p}d\mu (y)\right) ^{1/p}& \leq
\left(\dint_{B(x,t)}d(x,y)^{p}|%
\varrho (x)+\varrho (y)|^{p}d\mu (y)\right) ^{1/p} \\
& \leq \left(\dint_{B(x,t)}d(x,y)^{p}%
\left( \varrho (x)^{p}+\varrho (y)^{p}\right) d\mu (y)\right) ^{1/p} \\
& \leq t\left( \varrho (x)+T_{t}^{p}\varrho (x)\right) ,
\end{align*}
hence
\begin{eqnarray*}
\mathcal{E}_{X^{(p)}}(h,t) &\leq &t\left\| \varrho \right\|
_{X^{(p)}}+t\left\| T_{t}^{p}\varrho \text{\ }\right\| _{X^{(p)}} \\
&\preceq &t\left\| \varrho \right\| _{X^{(p)}}\text{ \ (by
(\ref{asi})).}
\end{eqnarray*}
In consequence
\begin{equation*}
\mathcal{E}_{X^{(p)}}(h,t)\preceq t\inf_{\varrho \in D(u)}\left\|
\varrho \right\| _{X^{(p)}}=t\left\| h\right\|
_{\dot{M}^{1,X^{(p)}}\left( \Omega \right) },
\end{equation*}
and taking the infimum over all representations of $f$ in $X^{(p)}+\dot{M}%
^{1,X^{(p)}}\left( \Omega \right) $, we obtain
\begin{equation*}
\mathcal{E}_{X^{(p)}}(f,t)\preceq
K(f,t,X^{(p)},\dot{M}^{1,X^{(p)}}).
\end{equation*}
\end{proof}

Having determined the $K-$functional between $X$ and \r{M}$^{1,X}$ in terms of
the $X-$modulus of smoothness, it is a known routine to determine the
corresponding interpolation spaces.

\begin{corollary}
Suppose $X$ is an r.i. space on $\Omega$ and $p\geq1.$ If $0<\theta<1$, $E$ is
a parameter and $b\in SV$, then
\begin{equation}
\text{\r{B}}_{X^{(p)},E}^{\theta,b}=(X^{(p)},\text{\r{M}}^{1,X^{(p)}}%
)_{\theta,b,E} \label{interhomo}%
\end{equation}
and
\begin{equation}
B_{X^{(p)},E}^{\theta,b}=(X^{(p)},M^{1,X^{(p)}})_{\theta,b,E},
\label{internohomo}%
\end{equation}
with equivalent norms.
\end{corollary}

\begin{proof}
Let us write $K(f,t,X^{(p)},$\r{M}$^{1,X^{(p)}})=K(f,t)$. Obviously
\[
\left\|  f\right\|  _{\text{\r{B}}_{X^{(p)},E}^{\theta,b}}=\left\|
t^{-\theta}b(t)E_{X^{(p)}}(f,t)\right\|  _{\tilde{E}}\preceq\left\|
t^{-\theta}b(t)K(f,t)\right\|  _{\tilde{E}}\leq\left\|  f\right\|
_{(X^{(p)},\text{\r{M}}^{1,X^{(p)}})_{\theta,b,E}}.
\]
Conversely,
\begin{align*}
\left\|  f\right\|  _{(X^{(p)},\text{\r{M}}^{1,X^{(p)}})_{\theta,b,E}}  &
=\left\|  t^{-\theta}b(t)K(f,t)\right\|  _{\tilde{E}}\preceq\left\|
t^{-\theta}b(t)\sum_{j=0}^{\infty}2^{-j}E_{X^{(p)}}(f,2^{j}t)\right\|
_{\tilde{E}}\\
&  \leq\sum_{j=0}^{\infty}2^{j(\theta-1)}\sup_{t>0}\frac{b(t)}{b(2^{j}%
t)}\left\|  \left(  2^{j}t\right)  ^{-\theta}b(2^{j}t)E_{X^{(p)}}%
(f,2^{j}t)\right\|  _{\tilde{E}}\\
&  \preceq\left(  \sum_{j=0}^{\infty}2^{j(\theta-1)}\sup_{t>0}\frac
{b(t)}{b(2^{j}t)}\right)  \left\|  t^{-\theta}b(t)E_{X^{(p)}}(f,t)\right\|
_{\tilde{E}}\text{ (by (\ref{propi})).}%
\end{align*}
To complete the proof, we need to see that the above series is convergent.
Consider $\varepsilon>0$ such that $\theta-1+\varepsilon<0,$ then it follows
from (\ref{sv1}) that
\[
\sup_{t>0}\frac{b(t)}{b(2^{j}t)}\preceq2^{j\varepsilon}%
\]
which implies convergence of the series.

To see (\ref{internohomo}), we observe that from Lemma \ref{pamam}, $\left\|
t^{-\theta}b(t)\min(1,t)\right\|  _{\tilde{L}^{1}}$ and $\left\|  t^{-\theta
}b(t)\min(1,t)\right\|  _{L^{\infty}}\ $ are both finite, so by interpolation
$\left\|  t^{-\theta}b(t)\min(1,t)\right\|  _{\tilde{E}}$ $<\infty,$ hence
\[
(X^{(p)},M^{1,X^{(p)}})_{\theta,b,E}\simeq(X^{(p)},\text{\r{M}}^{1,X^{(p)}%
})_{\theta,b,E}+\left\|  f\right\|  _{X^{(p)}}.
\]
Therefore
\begin{align*}
\left\|  f\right\|  _{B_{X^{(p)},E}^{\theta,b}}  &  =\left\|  t^{-\theta
}b(t)E_{X^{(p)}}(f,t)\right\|  _{\tilde{E}(0,1)}+\left\|  f\right\|
_{X^{(p)}}\\
&  \leq\left\|  t^{-\theta}b(t)K(f,t)\right\|  _{\tilde{E}}+\left\|
f\right\|  _{X^{(p)}}\\
&  \simeq\left\|  f\right\|  _{(X^{(p)},M^{1,X^{(p)}})_{\theta,b,E}}.
\end{align*}
Conversely, $E_{X^{(p)}}(f,t)\leq\left\|  f\right\|  _{X^{(p)}}$ gives
\begin{align*}
\left\|  t^{-\theta}b(t)K(f,t)\right\|  _{\tilde{E}}  &  \simeq\left\|
t^{-\theta}b(t)E_{X^{(p)}}(f,t)\right\|  _{\tilde{E}}\\
&  \simeq\left\|  t^{-\theta}b(t)E_{X^{(p)}}(f,t)\right\|  _{\tilde{E}%
(0,1)}+\left\|  t^{-\theta}b(t)E_{X^{(p)}}(f,t)\right\|  _{\tilde{E}%
(1,\infty)}\\
&  \preceq\left\|  t^{-\theta}b(t)E_{X^{(p)}}(f,t)\right\|  _{\tilde{E}%
(0,1)}+\left\|  f\right\|  _{X^{(p)}}\left\|  t^{-\theta}b(t)\right\|
_{\tilde{E}(1,\infty)}\\
&  \preceq\left\|  t^{-\theta}b(t)E_{X^{(p)}}(f,t)\right\|  _{\tilde{E}%
(0,1)}+\left\|  f\right\|  _{X^{(p)}}\ \text{(by Lemma \ref{pamam})}\\
&  =\left\|  f\right\|  _{B_{X^{(p)},E}^{\theta,b}}.
\end{align*}
\end{proof}

\begin{remark}
\label{equal} It is clear, that if in the previous corollary we use
the functional $\mathcal{E}_{X^{(p)}}(f,t)$ instead of
$E_{X^{(p)}}(f,t),$ then the same conclusion holds. In particular,
since for $1<p<\infty,$ $L^{p}=\left(  L^{1}\right)  ^{(p)},$ we get
\[
\text{\r{B}}_{L^{p},q}^{\theta}=\left(  L^{p},\dot{M}^{1,p}\right)
_{\theta,q}=\left(  \left(  L^{1}\right)  ^{(p)},\dot{M}^{1,p}\right)
_{\theta,q}=\mathcal{\dot{B}}_{p,q}^{s},
\]
with equivalent norms. Similarly $B_{L^{p},q}^{\theta}=\mathcal{B}_{p,q}^{s}. $
\end{remark}

A direct application of Theorem 5.1 and Corollary 5.2 of
\cite{FerSig} gives us

\begin{theorem}
Suppose $X$ is an r.i. space on $\Omega$, $E_{0},E_{1}$ and $E$
parameters, and $b,b_{0}$ and $b_{1}$ $SV-$ functions. Then, for
$0<\theta_{0}<\theta_{1}<1$ and $0<\theta<1$, we get
\[
\left(  \text{\r{B}}_{X,E_{0}}^{\theta_{0},b_{0}},\text{\r{B}}_{X,E_{1}%
}^{\theta_{1},b_{1}}\right)  _{\theta,b,E}=\text{\r{B}}_{X,E}^{\tilde{\theta
},\tilde{b}}%
\]
where $\tilde{\theta}=(1-\theta)\theta_{0}+\theta\theta_{1}$ and $\tilde
{b}(t)=b_{0}^{1-\theta}(t)b_{1}^{\theta}(t)b\left(  t^{\theta_{1}-\theta_{0}%
}b_{0}(t)/b_{1}(t)\right)  ,$ and
\[
\left(  \text{\r{M}}^{1,X},\text{\r{B}}_{X,E_{1}}^{\theta_{1},b_{1}}\right)
_{\theta,b,E}=\text{\r{B}}_{X,E}^{\hat{\theta},\hat{b}}%
\]
where $\hat{\theta}=\theta\theta_{1}$ and $\hat{b}(t)=b_{1}^{\theta
}(t)b\left(  t^{\theta_{1}}/b_{1}(t)\right)  .$

A similar result holds in the inhomogeneous case.
\end{theorem}

\subsection{Sobolev embedding theorems for Haj{\l}asz-Besov spaces.}\label{subsection41}

Throughout what follows we will write $K(f,t)$ instead of $K(f,t,X,$\r
{M}$^{1,X})$.

\begin{theorem}
\label{the1}Let $\left(  \Omega,d,\mu\right)  $ be an $RD-$space with indices
$\left(  k,n\right)  $ that satisfy the non-collapsing condition. Let $R(t)=\max(t^{1/n},t^{1/k})\ $and $f\in$\r{M}%
$^{1,L^{1}+L^{\infty}}$. Then for all $t>0,$ we have
\begin{equation}
O(f,t)\preceq\frac{K(f,R(t))}{\phi_{X}(t)}. \label{jaja2}%
\end{equation}
\end{theorem}

\begin{proof}
Fix $t>0$ and a decomposition $f=g+h$ with $g\in X$ and
$h\in\dot{M}^{1,X}$, since $f^{\ast}(t)\leq
g^{\ast}(t/2)+h^{\ast}(t/2)$ and $f^{\ast}(t)\geq h^{\ast
}(2t)-g^{\ast}(t)$ we have that
\begin{align}
f^{\ast\ast}\left(  t\right)  -f^{\ast}\left(  t\right)   &  \leq g^{\ast\ast
}\left(  t\right)  +h^{\ast\ast}\left(  t\right)  +g^{\ast}(t)-h^{\ast
}(2t)\label{fin}\\
&  \leq\left(  g^{\ast\ast}\left(  t\right)  +g^{\ast}(t)\right)  +\left(
h^{\ast\ast}\left(  t\right)  -h^{\ast\ast}\left(  2t\right)  \right)
+\left(  h^{\ast\ast}\left(  2t\right)  -h^{\ast}(2t)\right) \nonumber\\
&  \leq2g^{\ast\ast}\left(  t\right)  +2\left(  h^{\ast\ast}\left(  2t\right)
-h^{\ast}(2t)\right)  \text{.}\nonumber
\end{align}
On the other hand, by H\"{o}lder's inequality
\begin{equation}
g^{\ast\ast}\left(  t\right)  =\frac{1}{t}\int_{0}^{\infty}g^{\ast}\left(
s\right)  \chi_{\lbrack0,t]}(s)ds\leq\left\|  g\right\|  _{X}\frac
{\phi_{X^{\prime}}(t)}{t}=\frac{\left\|  g\right\|  _{X}}{\phi_{X}(t)}.
\label{f1}%
\end{equation}
By \cite[Theorem 1]{MarOrt}, it follows that for all $\rho\in D(h)\ $and for
all $t>0$, we have
\begin{equation*}
h^{\ast\ast}\left(  2t\right)  -h^{\ast}(2t)\preceq R(t)g^{\ast\ast}(t).
\end{equation*}
So, using the same argument as in (\ref{f1})
\begin{align}
h^{\ast\ast}\left(  2t\right)  -h^{\ast}(2t)  &  \preceq R(2t)\rho^{\ast\ast
}(2t)\label{f2}\\
&  \preceq R(t)\rho^{\ast\ast}(t)\nonumber\\
&  \leq R(t)\frac{\left\|  \rho\right\|  _{X}}{\phi_{X}(t)}\nonumber\\
&  =R(t)\frac{\left\|  h\right\|  _{\text{\r{M}}^{1,X}}}{\phi_{X}%
(t)}.\nonumber
\end{align}
Putting (\ref{f2}) and (\ref{f1}) back into (\ref{fin}), we find that
\[
f^{\ast\ast}\left(  t\right)  -f^{\ast}\left(  t\right)  \preceq\frac{1}%
{\phi_{X}(t)}\left(  \left\|  g\right\|  _{X}+R(t)\left\|  h\right\|
_{\dot{M}^{1,X}}\right)  .
\]
Taking the infimum over all decompositions of $f=g+h$ as shown above gives the
required estimate (\ref{jaja2}).
\end{proof}

The Theorem \ref{the1} is the key to prove the following
Sobolev-type embedding results for generalised Haj\l asz-Besov
spaces.

\begin{theorem}
\label{inclusi} Let $\left( \Omega ,d,\mu \right) $ be an $RD-$space
with indices $\left( k,n\right) $ which satisfies the non-collapsing
condition. Let $E$ be a parameter, $
b\in SV$, $0<\theta <1$ and $X$ be an r.i. space on $\Omega $. Then for all $f\in L^{1}+L^{\infty }$ such that%
\footnote{%
Condition $f^{\ast }(\infty )=0$ is equivalent to $\mu
\{|f|>t\}<\infty $ for all $t>0$.} $f^{\ast }(\infty )=0$, we have
\begin{equation*}
\left\| \frac{b(t^{1/n})\phi _{X}(t)}{t^{\theta /n}}O(f,t)\right\| _{\tilde{E%
}(0,1)}+\left\| \frac{b(t^{1/k})\phi _{X}(t)}{t^{\theta /k}}O(f,t)\right\| _{%
\tilde{E}(1,\infty )}\preceq \left\| f\right\|
_{\text{\r{B}}_{X,E}^{\theta ,b}}.
\end{equation*}
Moreover:

\begin{enumerate}
\item  If $\underline{\beta }_{X}>\frac{\theta }{k},$ then
\begin{equation}
\left\| \frac{b(t^{1/n})\phi _{X}(t)}{t^{\theta /n}}f^{\ast \ast
}\left(
t\right) \right\| _{\tilde{E}(0,1)}+\left\| \frac{b(t^{1/k})\phi _{X}(t)}{%
t^{\theta /k}}f^{\ast \ast }\left( t\right) \right\|
_{\tilde{E}(1,\infty )}\preceq \left\| f\right\|
_{\text{\r{B}}_{X,E}^{\theta ,b}}.  \label{todo}
\end{equation}

\item  If $\frac{\theta}{n}<\underline{\beta}_{X}\leq\frac{\theta}{k},$ then
\begin{equation*}
\left\| \frac{b(t^{1/n})\phi_{X}(t)}{t^{\theta/n}}f^{\ast\ast}\left(
t\right) \right\| _{\tilde{E}(0,1)}+\left\| \frac{b(t^{1/k})\phi_{X}(t)}{%
t^{\theta/k}}O(f,t)\right\| _{\tilde{E}(1,\infty)}\preceq\left\| f\right\| _{%
\text{\r{B}}_{X,E}^{\theta,b}}+\left\| f\right\|
_{L^{1}+L^{\infty}}.
\end{equation*}

\item  If $\underline{\beta }_{X}=\frac{\theta }{n},$ then

\begin{enumerate}
\item  If $\left\| \frac{t^{\theta /n}}{b(t^{1/n})\phi _{X}(t)}\right\| _{%
\tilde{E}^{\prime }(0,1)}<\infty ,$ then
\begin{equation*}
\left\| f\right\| _{L^{\infty }}\preceq \left\| f\right\| _{\text{\r{B}}%
_{X,E}^{\theta ,b}}+\left\| f\right\| _{L^{1}+L^{\infty }}.
\end{equation*}

\item  If the function $b(t^{1/n})\phi _{X}(t)/t^{\mathbb{\theta }/n}$ is
quasi-increasing and Boyd's indices \underline{$\alpha $}$_{E},\bar{\alpha}%
_{E}\in (0,1),$ then
\begin{equation*}
\left\| \frac{b(t^{1/n})\phi _{X}(t)}{t^{\theta /n}\ell (t)}f^{\ast
\ast
}\left( t\right) \right\| _{\tilde{E}(0,1)}\preceq \left\| f\right\| _{\text{%
\r{B}}_{X,E}^{\theta ,b}}+\left\| f\right\| _{L^{1}+L^{\infty }}.
\end{equation*}

\item  If \begin{equation*}\sup_{t\in (0,1)}\frac{t^{\theta /n}}{b(t^{1/n})\phi _{X}(t)}%
<\infty ,\end{equation*} then\footnote{%
The local spaces of functions of bounded mean oscillation, $bmo$,
consists of all locally integrable functions satisfying the
condition
\begin{equation*}
bmo(\Omega )=\left\{ f:\sup_{\mu (B)\leq 1}\dint_{B}\left|
f(y)-\dint_{B}f(s)d\mu (s)\right| d\mu (y)+\sup_{\mu
(B)>1}\dint_{B}f(s)d\mu (s)\right\}
\end{equation*}
where $B$ is a ball in $\Omega .$ This definition coincides with
\cite[2.2.2 (viii), p. 37]{Tri83}.}
\begin{equation*}
\left\| f\right\| _{bmo}\preceq \left\| f\right\| _{\text{\r{B}}%
_{X,L^{\infty }}^{\theta ,b}}+\left\| f\right\| _{L^{1}+L^{\infty
}}.
\end{equation*}
\end{enumerate}

\item  If $\bar{\beta}_{X}<\frac{\theta }{n},$ then
\begin{equation*}
\left\| f\right\| _{L^{\infty }}\preceq \left\| f\right\| _{\text{\r{B}}%
_{X,E}^{\theta ,b}}+\left\| f\right\| _{L^{1}+L^{\infty }}.
\end{equation*}
\end{enumerate}
\end{theorem}

\begin{proof}
It follows from (\ref{jaja2}) that
\begin{equation*}
O(f,t)\phi _{X}(t)\preceq K(f,R(t))
\end{equation*}
thus
\begin{equation*}
\frac{b(R(t))\phi _{X}(t)}{R(t)^{\theta }}O(f,t)\preceq b(R(t))\frac{%
K(f,R(t))}{R(t)^{\theta }}.
\end{equation*}
Whence
\begin{align}
\left\| \frac{b(t^{1/n})\phi _{X}(t)}{t^{\mathbb{\theta
}/n}}O(f,t)\right\|
_{\tilde{E}(0,1)}+\left\| \frac{b(t^{1/k})\phi _{X}(t)}{t^{\mathbb{\theta }%
/k}}O(f,t)\right\| _{\tilde{E}(1,\infty )}& \simeq \left\|
\frac{b(R(t))\phi
_{X}(t)}{R(t)^{\theta }}O(f,t)\right\| _{\tilde{E}}  \label{calcul} \\
& \preceq \left\| b(R(t))\frac{K(f,R(t))}{R(t)^{\theta }}\right\| _{\tilde{E}%
}  \notag \\
& \preceq \left\| b(t)\frac{K(f,t)}{t^{\theta }}\right\|
_{\tilde{E}}\text{
\ (by (\ref{propi}))}  \notag \\
& =\left\| f\right\| _{\text{\r{B}}_{X,E}^{\theta ,b}}.  \notag
\end{align}
and we have the inequality (\ref{todo}).

\textbf{(i)} Let $\varphi (t)=\frac{b(R(t))\phi _{X}(t)}{%
R(t)^{\theta }},$ by (\ref{calcul}) we know that
\begin{equation*}
\left\| \varphi (t)O(f,t)\right\| _{\tilde{E}}\preceq \left\| f\right\| _{%
\text{\r{B}}_{X,E}^{\theta ,b}}.
\end{equation*}
Taking into account (\ref{der2est}) and $f^{\ast }\left( \infty
\right)=0$, from the fundamental theorem of calculus we deduce that
\begin{equation*}
f^{\ast \ast }(t)=\int_{t}^{\infty }O(f,s)\frac{ds}{s},
\end{equation*}
thus we need to see that
\begin{equation}
\left\| \varphi (t)f^{\ast \ast }(t)\right\| _{\tilde{E}}=\left\|
\varphi (t)\int_{t}^{\infty }O(f,s)\frac{ds}{s}\right\|
_{\tilde{E}}\preceq \left\| \varphi (t)O(f,t)\right\| _{\tilde{E}}.
\label{cota1}
\end{equation}
To do this, we look at the operator
\begin{equation*}
(Ug)(t)=\varphi (t)\int_{t}^{\infty }g(z)\frac{1}{\varphi (z)}\frac{dz}{z},%
\text{ \ }(t>0),
\end{equation*}
and we will prove that it is bounded for any parameter $\tilde{E}$.
Due to the interpolation properties of the space $\tilde{E}$, it
suffices to prove that $U$ is bounded in $\tilde{L}^{1}$ and in
$L^{\infty }.$

Since $b(R(t))\in SV,$ it follows from the proposition \ref{sv} that $%
\underline{\beta}_{\varphi}=-\frac{\theta}{k}+\underline{\beta}_{X}>0$ and $%
\overline{\beta}_{1/\varphi}=-\underline{\beta}_{\varphi}<0$, hence
by (\ref {indi})
\begin{equation*}
\int_{0}^{t}\varphi(z)\frac{dz}{z}\preceq\varphi(t)\text{ and }\int
_{t}^{\infty}\frac{1}{\varphi(z)}\frac{dz}{z}\preceq\frac{1}{\varphi(t)}.
\end{equation*}
Consequently
\begin{align*}
\left\| Ug\right\| _{\tilde{L}^{1}} & \leq\int_{0}^{\infty}\left(
\varphi(t)\int_{t}^{\infty}\left| g(z)\right|
\frac{1}{\varphi(z)}\frac
{dz}{z}\right) \frac{dt}{t} \\
& =\int_{0}^{\infty}\left| g(t)\right| \varphi(t)\left( \int_{0}^{t}\frac{1}{%
\varphi(z)}\frac{dz}{z}\right) \frac{dt}{t}\preceq\left\| g\right\| _{\tilde{%
L}^{1}};
\end{align*}
and
\begin{equation*}
\left\| Ug\right\| _{L^{\infty}}=\sup_{t>0}\int_{t}^{\infty}\left|
g(z)\right| \varphi(t)\frac{dz}{z}\leq\left\| Ug\right\|
_{L^{\infty}}\sup_{t>0}\int_{t}^{\infty}\frac{1}{\varphi(z)}\frac{dz}{z}%
\preceq\left\| g\right\| _{L^{\infty}}.
\end{equation*}
If we take $g(t)=O(f,t)\varphi(t)$, we have the desired inequality
(\ref {cota1}).

\textbf{(ii)} Using an argument analogous to the previous case, it
follows that the operator
\begin{equation*}
(Ug)(t)=\frac{b(t^{1/n})\phi _{X}(t)}{t^{\mathbb{\theta }/n}}\int_{t}^{1}g(z)%
\frac{z^{\mathbb{\theta }/n}}{b(z^{1/n})\phi _{X}(z)}\frac{dz}{z}
\end{equation*}
is bounded in $\tilde{E}(0,1).$

Taking $g(t)=\frac{b(t^{1/n})\phi _{X}(t)}{t^{\theta /n}}\left(
f^{\ast \ast }\left( t\right) -f^{\ast }\left( t\right) \right) ,$
we have that
\begin{align}
Uh(t)& =\frac{b(t^{1/n})\phi _{X}(t)}{t^{\theta
/n}}\int_{t}^{1}\left( f^{\ast \ast }\left( z\right) -f^{\ast
}\left( z\right) \right) \frac{dz}{z}
\label{esto} \\
& =\frac{b(t^{1/n})\phi _{X}(t)}{t^{\theta /n}}\left( f^{\ast \ast
}\left(
t\right) -f^{\ast \ast }(1)\right)  \notag \\
& =\frac{b(t^{1/n})\phi _{X}(t)}{t^{\theta /n}}\left( f^{\ast \ast
}\left( t\right) -\left\| f\right\| _{L^{1}+L^{\infty }}\right)
\notag
\end{align}
thus
\begin{align*}
\left\| \frac{b(t^{1/n})\phi _{X}(t)}{t^{\theta /n}}f^{\ast \ast
}\left(
t\right) \right\| _{\tilde{E}(0,1)}& \preceq \left\| Uh\right\| _{\tilde{E}%
(0,1)} \\
& \preceq \left\| \frac{b(t^{1/n})\phi _{X}(t)}{t^{\theta
/n}}O(f,t)\right\|
_{\tilde{E}(0,1)}+\left\| f\right\| _{L^{1}+L^{\infty }} \\
& \preceq \left\| f\right\| _{\text{\r{B}}_{X,E}^{\theta
,b}}+\left\| f\right\| _{L^{1}+L^{\infty }}.
\end{align*}

\textbf{(iii) case (a) } From the fundamental theorem of calculus
and the H\"{o}lder inequality, we obtain
\begin{align*}
f^{\ast \ast }(t)-\left\| f\right\| _{L^{1}+L^{\infty }}& =\int_{t}^{1}O(f,t)%
\frac{b(t^{1/n})\phi _{X}(t)}{t^{\mathbb{\theta
}/n}}\frac{t^{\mathbb{\theta
}/n}}{b(t^{1/n})\phi _{X}(t)}\frac{dt}{t} \\
& \leq \left\| O(f,t)\frac{b(t^{1/n})\phi _{X}(t)}{t^{\mathbb{\theta }%
/n}}\right\| _{\tilde{E}(0,1)}\left\| \frac{s^{\mathbb{\theta }/n}}{%
b(s^{1/n})\phi _{X}(s)}\chi _{\lbrack t,1]}(s)\right\|
_{\tilde{E}^{^{\prime }}(0,1)}\text{(by (\ref{Hold}))}
\end{align*}
thus
\begin{equation*}
\left\| f\right\| _{L^{\infty }}=\sup_{0<t<1}f^{\ast \ast
}(t)\preceq \left\| f\right\| _{\text{\r{B}}_{X,E}^{\theta
,b}}+\left\| f\right\| _{L^{1}+L^{\infty }},
\end{equation*}
as we want to show.

\textbf{(iii) case (b) }Consider the operator
\begin{equation*}
\left( Ug\right) (t)=\frac{b(t^{1/n})\phi _{X}(t)}{t^{\mathbb{\theta }%
/n}\ell (t)}\int_{t}^{1}g(z)\frac{dz}{z}.
\end{equation*}
If $1<p<\infty ,$ then by \cite[Theorem 6.5]{BR} we have
\begin{align*}
\left\| Ug\right\| _{\tilde{L}^{p}(0,1)}& \preceq \left(
\int_{0}^{1}\left(
\frac{1}{\ell (t)}\int_{t}^{1}\left| g(z)\right| \frac{b(z^{1/n})\phi _{X}(z)%
}{z^{\mathbb{\theta }/n}}\frac{dz}{z}\right) ^{p}\frac{dt}{t}\right)
^{1/p}
\\
& \leq \left( \int_{0}^{1}\left( g(t)\frac{b(t^{1/n})\phi _{X}(t)}{t^{%
\mathbb{\theta }/n}}\right) ^{p}\frac{dt}{t}\right) ^{1/p} \\
& =\left\| g(t)\frac{b(t^{1/n})\phi _{X}(t)}{t^{\mathbb{\theta
}/n}}\right\| _{\tilde{L}^{p}(0,1)},
\end{align*}
and
\begin{align*}
\left\| Ug\right\| _{L^{\infty }[0,1]}& \preceq \sup_{0<t<1}\left( \frac{1}{%
\ell (t)}\int_{t}^{1}\left| g(z)\right| \frac{b(z^{1/n})\phi _{X}(z)}{z^{%
\mathbb{\theta }/n}}\frac{dz}{z}\right) \\
& \leq \sup_{0<t<1}\left( \sup_{t<s<1}\left| g(t)\frac{b(t^{1/n})\phi _{X}(t)%
}{t^{\mathbb{\theta }/n}}\right| \left( \frac{1}{\ell (t)}\int_{t}^{1}\frac{%
dz}{z}\right) \right) \\
& \preceq \left\| g(t)\frac{b(t^{1/n})\phi _{X}(t)}{t^{\mathbb{\theta }/n}}%
\right\| _{L^{\infty }[0,1]}.
\end{align*}
Choose $1<p<\infty $ such $\overline{\alpha }_{E}<1/p$, since also
\underline{$\alpha $}$_{E}>0$, $E$ is an interpolation space between $%
L^{p}[0,1]$ and $L^{\infty }[0,1]$ (see \cite[Theorem 2.b.11]{LT}),
therefore
\begin{equation*}
\left\| Ug\right\| _{\tilde{E}(0,1)}\preceq \left\|
\frac{b(t^{1/n})\phi _{X}(t)}{t^{\mathbb{\theta }/n}}g(t)\right\|
_{\tilde{E}(0,1)}.
\end{equation*}
Taking $g(t)=O(f,t)$, we get
\begin{align*}
\left\| \frac{b(t^{1/n})\phi _{X}(t)}{\ell (t)t^{\theta /n}}f^{\ast
\ast
}\left( t\right) \right\| _{\tilde{E}(0,1)}& \preceq \left\| \frac{%
b(t^{1/n})\phi _{X}(t)}{t^{\mathbb{\theta }/n}}O(f,t)\right\| _{\tilde{E}%
(0,1)}+\left\| f\right\| _{L^{1}+L^{\infty }} \\
& \preceq \left\| f\right\| _{\text{\r{B}}_{X,E}^{\theta
,b}}+\left\| f\right\| _{L^{1}+L^{\infty }}.
\end{align*}

\textbf{(iii) case (c) }Given $B=B(x,r)$ a ball centered at $x$ with radius $%
r, $ we have that
\begin{align*}
I(B)& :=\dint_{B}\left| f(y)-
\dint_{B}f(s)d\mu (s)\right| d\mu (y) \\
& \leq \dint_{B}\left( \dint_{B}\left|
f(y)-f(s)\right| d\mu (s)\right) d\mu (y) \\
& \leq \frac{E_{X}(f,r)}{\mu (B)}\left\| \chi _{B}\right\| _{X^{\prime }}%
\text{ (by (\ref{Hold}))} \\
& =\frac{E_{X}(f,r)}{\phi _{X}(\mu (B))} \\
& \preceq \frac{K(f,r)}{\phi _{X}(\mu (B))}\text{ (by Theorem \ref
{TeoInterpol})}.
\end{align*}
We claim that given $t>0$ and $x\in \Omega ,$ there exists a ball
$B(x)$ centered at $x$ such that
\begin{equation}
t/C_{D}\leq \mu (B(x))\leq t  \label{aa}
\end{equation}
furthermore, for any ball $B(x,r)$ satisfying (\ref{aa}), there
exists a constant $c=c(c_{0},\kappa ,n,k),$ such that
\begin{equation}
r\leq cR(t).  \label{claim}
\end{equation}
Assuming for the moment the validity of (\ref{aa}) and
(\ref{claim}), and considering that that $K(f,\cdot )$ is
increasing, we get
\begin{align}
\sup_{\mu (B)\leq 1}I(B)& \preceq \sup_{\mu (B)\leq
1}\frac{K(f,r)}{\phi
_{X}(\mu (B))}\leq \sup_{0<t\leq 1}\sup_{t/C_{D}\leq \mu (B)\leq t\leq 1}%
\frac{K(f,r)}{\phi _{X}(\mu (B))}  \label{aq} \\
& \leq \sup_{0<t\leq 1}\sup_{t/C_{D}\leq \mu (B)\leq t\leq 1}\frac{%
K(f,ct^{1/n})}{\phi _{X}(\mu (B))}\leq \sup_{0<t\leq 1}\frac{K(f,ct^{1/n})}{%
\phi _{X}(t/C_{D})}  \notag \\
& \preceq \sup_{0<t\leq 1}\frac{K(f,t^{1/n})}{\phi _{X}(t)}  \notag \\
& \leq \sup_{0<t\leq 1}\frac{K(f,t^{1/n})b(t^{1/n})}{t^{\theta
/n}}\left( \sup_{0<t\leq 1}\frac{t^{\theta /n}}{\phi
_{X}(t)b(t^{1/n})}\right)  \notag
\\
& \preceq \left\| f\right\| _{\text{\r{B}}_{X,L^{\infty }}^{\theta
,b}}, \notag
\end{align}
where in the third line of the above inequality we have used the
concavity of the functions $K(f,\cdot )$ and $\phi _{X}(\cdot ).$

It remains to prove (\ref{aa}) and (\ref{claim}). Consider
$r_{0}=\sup
\left\{ r:\mu \left( B(x,r)\right) <t/C_{D}\right\} ,$ then $\ $%
\begin{equation*}
\mu \left( B(x,r_{0})\right) \leq t/C_{D}\leq \mu \left(
B(x,2r_{0})\right) \leq C_{D}\mu \left( B(x,r_{0})\right) \leq t.
\end{equation*}
If $B(x,r)$ satisfies (\ref{aa}), then obviously there exists a constant $%
c=c(c_{0},\kappa ,n,k)$ such that $R(t/c_{0}\kappa )\leq cR(t),$ so
it will be enough to see that
\begin{equation*}
r\leq R(t/c_{0}\kappa ).
\end{equation*}
If $\frac{t}{c_{0}\kappa }\leq 1,$ then $R(t/c_{0}\kappa )=\left(
t/c_{0}\kappa \right) ^{1/n},$ therefore, if $r>\left( t/c_{0}\kappa
\right) ^{1/n},$ then from (\ref{belw}) we will get $\mu
(B(x,r))>t,$ in contradiction to (\ref{aa}). In the case of
$t/c_{0}\kappa \geq 1$, we proceed in a similar way.

Finally
\begin{eqnarray*}
\left\| f\right\| _{bmo} &=&\sup_{\mu (B)\leq 1}I(B)+\sup_{\mu (B)\geq 1}%
\dint_{B}\left| f(x)\right| d\mu (x) \\
&\preceq &\left\| f\right\| _{\text{\r{B}}_{X,L^{\infty }}^{\theta
,b}}+\sup_{\mu (B)\geq 1}\dint_{B}\left| f(x)\right| d\mu (x)%
\text{ (by (\ref{aq}))} \\
&\leq &\left\| f\right\| _{\text{\r{B}}_{X,L^{\infty }}^{\theta
,b}}+\sup_{\mu (B)\geq 1}\frac{1}{\mu (B)}\int_{0}^{\mu (B)}f^{\ast }(s)ds \\
&=&\left\| f\right\| _{\text{\r{B}}_{X,L^{\infty }}^{\theta
,b}}+f^{\ast
\ast }(1) \\
&=&\left\| f\right\| _{\text{\r{B}}_{X,L^{\infty }}^{\theta
,b}}+\left\| f\right\| _{L^{1}+L^{\infty }}.
\end{eqnarray*}

\textbf{(iv) }Let $\varphi (t)=\frac{b(t^{1/n})\phi
_{X}(t)}{t^{\theta /n}}.$
Since $b(t^{1/n})\in SV$, from the Proposition \ref{sv} we get $\overline{%
\beta }_{\varphi }=-\frac{\theta }{n}+\overline{\beta }_{X}<0$ and $%
\underline{\beta }_{1/\varphi }=-\overline{\beta }_{\varphi }>0,$
therefore by (\ref{indi})
\begin{equation*}
\left\| \frac{1}{\varphi }\right\| _{\tilde{L}^{1}(0,t)}=\int_{0}^{t}\frac{1%
}{\varphi (s)}\frac{ds}{s}\preceq \frac{1}{\varphi (t)}.
\end{equation*}
With the choice of $\underline{\beta }_{1/\varphi }>a>0$ we have the
function $\frac{1}{\varphi (t)t^{a}}$ is almost increasing,
therefore
\begin{equation*}
\left\| \frac{1}{\varphi }\right\| _{L^{\infty }(0,t)}=\sup_{0<s<t}\frac{%
s^{a}}{\varphi (s)s^{a}}\preceq \frac{1}{\varphi (t)}.
\end{equation*}
Using interpolation, for any parameter $E$ we have
\begin{equation*}
\left\| \frac{1}{\varphi }\right\| _{\tilde{E}(0,t)}\preceq
\frac{1}{\varphi (t)},
\end{equation*}
and, by H\"{o}lder inequality
\begin{align*}
\left\| f\right\| _{L^{\infty }}& =f^{\ast \ast }(0)=\int_{0}^{1}O(f,t)\frac{%
dt}{t}+\left\| f\right\| _{L^{1}+L^{\infty }} \\
& \leq \left\| O(f,t)\varphi (t)\right\| _{\tilde{E}(0,1)}\left\| \frac{1}{%
\varphi (t)}\right\| _{\tilde{E}^{^{\prime }}(0,1)}+\left\|
f\right\|
_{L^{1}+L^{\infty }} \\
& \preceq \left\| f\right\| _{\text{\r{B}}_{X,E}^{\theta
,b}}+\left\| f\right\| _{L^{1}+L^{\infty }}
\end{align*}
as we wanted to see.
\end{proof}

Now we consider the essential continuity problem and obtain Morrey
type results for functions in generalised Haj\l asz-Besov spaces.

\begin{theorem}
\label{contt}Let $\left( \Omega ,d,\mu \right) $ be an $RD-$space
with indices $\left( k,n\right) $ which satisfies the non-collapsing
condition and $X$ an r.i. space on $\Omega $ and $f\in X\cap
L^{\infty }.$ Then
\begin{equation*}
\left| f(x)-f(y)\right| \preceq \int_{0}^{8d(x,y)}\frac{{E}_{X}(f,s)%
}{\varphi _{X}\left( \min (s^{k},s^{n})\right) }\frac{ds}{s}.
\end{equation*}
In particular, if
\begin{equation*}
\int_{0}^{1}\frac{{E}_{X}(f,s)}{\varphi _{X}\left( s^{n}\right) }%
\frac{ds}{s}<\infty ,
\end{equation*}
then, $f$ is essentially continuous.
\end{theorem}

\begin{proof}
Given $r>0$ and $f\in X\cap L^{\infty },$ let
\begin{equation*}
\nabla _{r}f(x):=\dint_{B(x,r)}\left| f(x)-f(y)\right| dy,
\end{equation*}
Since $\nabla _{r}f\in X\cap L^{\infty }$, it follows that
\begin{equation*}
\lim_{q\rightarrow \infty }E_{X^{(q)}}(f,r)=\lim_{q\rightarrow
\infty }\left\| \left( \nabla _{r}f\right) ^{q}\right\|
_{X}^{1/q}=\sup_{x\in \Omega }\nabla _{r}f(x).
\end{equation*}
On the other hand, from Lemma \ref{Suerte} (ii) we get (notice that
$\mathcal{E}_{X^{(1)}}(f,s)=E_X(f,s)$)
\begin{equation}
\lim_{q\rightarrow \infty }E_{X^{(q)}}(f,r)\leq c\int_{0}^{4r}\frac{{%
E}_{X}(f,s)}{\varphi _{X}\left( \min (s^{k},s^{n})\right)
}\frac{ds}{s}. \label{suerte 2}
\end{equation}
Let $B(x,r)$ be a ball centered at $x$ with radius $r$, since if
$y\in B(x,r) $, then $B(x,r)\subset B(y,2r)$\ we have
\begin{eqnarray*}
\left| f(x)-f(y)\right| &\leq &\dint_{B(x,r)}\left| f(x)-f(z)\right|
d\mu
(z)+\dint_{B(x,r)}\left| f(z)-f(y)\right| d\mu (z) \\
&\preceq &\dint_{B(x,2r)}\left| f(x)-f(z)\right| d\mu
(z)+\dint_{B(y,2r)}\left| f(z)-f(y)\right| d\mu (z)\text{ (since
}\mu \text{
is doubling)} \\
&\leq &\sup_{x\in \Omega }\dint_{B(x,2r)}\left| f(x)-f(z)\right|
d\mu
(z)+\sup_{y\in \Omega }\dint_{B(y,2r)}\left| f(z)-f(y)\right| d\mu (z) \\
&\preceq &\int_{0}^{8r}\frac{{E}_{X}(f,s)}{\varphi _{X}\left( \min
(s^{k},s^{n})\right) }\frac{ds}{s}\text{ (by (\ref{suerte 2})).}
\end{eqnarray*}
So, if $x,y\in \Omega $ are such that $r/8=d(x,y)<1,$ then
\begin{equation*}
\left| f(x)-f(y)\right| \preceq \int_{0}^{d(x,y)}\frac{{E}_{X}(f,s)}{%
\varphi _{X}\left( s^{n}\right) }\frac{ds}{s}
\end{equation*}
and the essential continuity of $f$ follows.
\end{proof}

\begin{corollary}
Let $\left( \Omega ,d,\mu \right) $ be an $RD-$space with indices
$\left( k,n\right) $ which satisfies the non-collapsing condition.
Let $E$ be a parameter, $b\in SV$, $0<\theta <1$ and $X$ be an r.i.
space on $\Omega $. Let $f\in B_{X,E}^{\theta ,b},$ then

\begin{enumerate}
\item  If $\underline{\beta }_{X}=\frac{\theta }{n}$ and $\left\| \frac{%
t^{\theta /n}}{b(t^{1/n})\phi _{X}(t)}\right\| _{\tilde{E}^{\prime
}(0,1)}<\infty ,$ then $f$ is essentially continuous.

\item  If $\underline{\beta }_{X}<\frac{\theta }{n},$ then for any $\bar{%
\beta}_{X}\leq \gamma \leq 1$ such that the function $\frac{\phi _{X}(t)}{%
t^{\gamma }}$ is quasi decreasing, the following Morrey type
embedding result holds:
\begin{equation*}
\ \left| f(x)-f(y)\right| \preceq \left\| f\right\| _{B_{X,E}^{\theta ,b}}%
\frac{d(x,y)^{\theta -\gamma n}}{b(d(x,r))}.
\end{equation*}
\end{enumerate}
\end{corollary}

\begin{proof}
(i) Since $\underline{\beta}_{X}=\frac{\theta}{n}$, by Theorem
\ref{inclusi} (iii) we have that $f\in L^{\infty }\cap X$.
Furthermore, by H\"{o}lder's inequality, we get
\begin{eqnarray*}
\int_{0}^{1}\frac{{E}_{X}(f,s)}{\varphi _{X}\left( s^{n}\right) }%
\frac{ds}{s} &=&\int_{0}^{1}s^{-\theta }b(s)\mathcal{E}_{X}(f,s)\frac{%
s^{\theta }}{b(s)\varphi _{X}\left( s^{n}\right) }\frac{ds}{s} \\
&\leq &\left\| s^{-\theta }b(s){E}_{X}(f,s)\right\| _{\tilde{E}%
(0,1)}\left\| \frac{s^{\theta }}{b(s)\varphi _{X}(s^{n})}\right\| _{\tilde{E}%
^{\prime }(0,1)} \\
&\preceq &\left\| s^{-\theta }b(s){E}_{X}(f,s)\right\| _{\tilde{E}%
(0,1)}\left\| \frac{t^{\theta /n}}{b(t^{1/n})\phi _{X}(t)}\right\| _{\tilde{E%
}^{\prime }(0,1)} \\
&\preceq &\left\| f\right\| _{B_{X,E}^{\theta ,b}}
\end{eqnarray*}
and the essential continuity of $f$ follows from the Theorem
\ref{contt}.

(ii) Since $\underline{\beta }_{X}<\frac{\theta }{n},$ by Theorem
\ref {inclusi} (iv) $f\in X\cap L^{\infty }$ so by Theorem
\ref{contt}$.$
\begin{equation*}
\left| f(x)-f(y)\right| \preceq \int_{0}^{8r}\frac{{E}_{X}(f,s)}{%
\varphi _{X}\left( \min (s^{k},s^{n})\right) }\frac{ds}{s}.
\end{equation*}
Considering that $\gamma n-\theta <0$ and that the function
$\frac{\phi _{X}(t)}{t^{\gamma }}$ is quasi decreasing, we get the
following
\begin{eqnarray*}
r^{\gamma n-\theta }b(r)\left| f(x)-f(y)\right|  &\preceq &r^{\gamma
n-\theta }b(r)\int_{0}^{8r}\frac{{E}_{X}(f,s)}{\varphi _{X}\left(
\min (s^{k},s^{n})\right) }\frac{ds}{s} \\
&\leq &\sup_{r>0}\left( \left( 8r\right) ^{\gamma n-\theta }b\left(
8r\right) \int_{0}^{8r}\frac{{E}_{X}(f,s)}{\varphi _{X}\left( \min
(s^{k},s^{n})\right) }\frac{ds}{s}\right)  \\
&\preceq &\sup_{r>0}\left( r^{\gamma n-\theta }b(r)\frac{{E}_{X}(f,r)%
}{\varphi _{X}\left( \min (r^{k},r^{n})\right) }\right) \text{ (by
Lemma \ref
{pamam})} \\
&\leq &\sup_{0<r<1}b(r)t^{\theta -\gamma n}\left( \frac{r^{\gamma n}}{%
\varphi _{X}\left( r^{n}\right) }\right) \frac{{E}_{X}(f,r)}{%
r^{\gamma n}}+\sup_{r\geq 1}b(r)r^{\theta -\gamma n}\frac{{E}%
_{X}(f,r)}{\varphi _{X}\left( r^{k}\right) } \\
&\preceq &\sup_{0<r<1}b(r)t^{-\theta }{E}_{X}(f,t)+\sup_{r\geq
1}b(r)r^{\theta -\gamma n}\frac{{E}_{X}(f,r)}{\varphi _{X}\left(
r^{k}\right) } \\
&\leq &\sup_{0<r<1}b(r)t^{-\theta }{E}_{X}(f,t)+\sup_{r\geq
1}b(r)r^{\theta -\gamma n}{E}_{X}(f,r) \\
&\preceq &\left\| f\right\| _{\text{\r{B}}_{X,L^{\infty }}^{\theta
,b}}+\left\| f\right\| _{X}\text{ \ (by (\ref{ultim}))} \\
&=&\left\| f\right\| _{B_{X,L^{\infty }}^{\theta ,b}} \\
&\preceq &\left\| f\right\| _{B_{X,E}^{\theta ,b}}\text{ (by
Proposition \ref {emm}).}
\end{eqnarray*}
Finally, if $x,y\in \Omega $ are such that $r=d(x,y)$, then
\begin{equation*}
\ \left| f(x)-f(y)\right| \preceq \left\| f\right\| _{B_{X,E}^{\theta ,b}}%
\frac{d(x,y)^{\theta -\gamma n}}{b(d(x,r))}.
\end{equation*}
\end{proof}

\begin{remark}
For finite measure metric spaces with convex isoperimetric profile,
a related result was obtained in \cite[Chapter 4]{Mami1}.
\end{remark}

We conclude this section with an example and some comments.
\begin{example}
Let $\left( \Omega ,d,\mu \right) $ be an $RD-$space with indices
$\left( k,n\right) $ that satisfies the non-collapsing condition.
Let $X$ be an r.i.
space on $\Omega ,$ $E=L^{q}$ ($1\leq q\leq \infty )$, $b\in SV$ and $%
0<\theta <1$. Assume that $\varphi _{X}(t)=t^{1/p}\psi (t)$ with
$1\leq p<\infty $ and $\psi \in SV.$ Then

\begin{enumerate}
\item  If $\frac{1}{p}>\frac{\theta }{k},$ then
\begin{align*}
& \left( \int_{0}^{1}\left( b(t^{1/n})\psi (t)f^{\ast \ast }\left(
t\right) t^{\frac{1}{p}-\frac{\theta }{n}}\right)
^{q}\frac{dt}{t}\right) ^{1/q}+\left( \int_{1}^{\infty }\left(
b(t^{1/n})\psi (t)f^{\ast \ast
}\left( t\right) t^{\frac{1}{p}-\frac{\theta }{k}}\right) ^{q}\frac{dt}{t}%
\right) ^{1/q} \\
& \preceq \left\| f\right\| _{\text{\r{B}}_{X,{L}^{q}}^{b,\theta }}.
\end{align*}

\item  If $\frac{\theta }{n}<\frac{1}{p}\leq \frac{\theta }{k},$ then
\begin{align*}
& \left( \int_{0}^{1}\left( b(t^{1/n})\psi (t)f^{\ast \ast }\left(
t\right) t^{\frac{1}{p}-\frac{\theta }{n}}\right)
^{q}\frac{dt}{t}\right)
^{1/q}+\left( \int_{1}^{\infty }\left( b(t^{1/n})\psi (t)O(f,t)t^{\frac{1}{p}%
-\frac{\theta }{k}}\right) ^{q}\frac{dt}{t}\right) ^{1/q} \\
& \preceq \left\| f\right\| _{\text{\r{B}}_{X,{L}^{q}}^{b,\theta
}}+\left\| f\right\| _{L^{1}+L^{\infty }.}
\end{align*}

\item  If $\frac{\theta }{n}=\frac{1}{p},$ then:

\begin{enumerate}
\item  If $\sup_{t\in (0,1)}\frac{1}{b\left( s^{1/n}\right) \psi (s)}<\infty
,$ then
\begin{equation*}
\left\| f\right\| _{L^{\infty }}\preceq \left\| f\right\| _{\text{\r{B}}%
_{X,L^{1}}^{b,\theta }}+\left\| f\right\| _{L^{1}+L^{\infty }.}
\end{equation*}
and $f$ is and essentially continuous.

\item  If $b\left( s^{1/n}\right) \psi (s)$ is quasi-increasing and $%
1<q<\infty ,$ then
\begin{equation*}
\left( \int_{0}^{1}\left( \frac{b(t^{1/n})\psi (s)}{1+\ln \left( \frac{1}{s}%
\right) }f^{\ast \ast }\left( s\right) \right)
^{q}\frac{ds}{s}\right) ^{1/q}\preceq \left\| f\right\|
_{\text{\r{B}}_{X,L^{q}}^{b,\theta }}+\left\| f\right\|
_{L^{1}+L^{\infty }.}
\end{equation*}

\item  If $0<\frac{1}{p}<1,$ then
\begin{equation*}
\left\| f\right\| _{bmo}\leq \left\| f\right\|
_{\text{\r{B}}_{X,L^{\infty }}^{b,\theta }}+\left\| f\right\|
_{L^{1}+L^{\infty },}
\end{equation*}
\end{enumerate}

\item  If $\frac{1}{p}<\frac{\theta }{n},$ then
\begin{equation*}
\left\| f\right\| _{L^{\infty }}\preceq \left\| f\right\| _{\text{\r{B}}%
_{X,L^{q}}^{b,\theta }}+\left\| f\right\| _{L^{1}+L^{\infty }}.
\end{equation*}
Moreover
\begin{equation*}
\ \left| f(x)-f(y)\right| \preceq \left\| f\right\| _{B_{X,E}^{\theta ,b}}%
\frac{d(x,y)^{\theta -\frac{n}{p}}}{b(d(x,r))}
\end{equation*}
in case that $\psi (s)$ is quasi-increasing, otherwise, for any $\frac{1}{p}%
<\gamma \leq 1$
\begin{equation*}
\ \left| f(x)-f(y)\right| \preceq \left\| f\right\| _{B_{X,E}^{\theta ,b}}%
\frac{d(x,y)^{\theta -\gamma n}}{b(d(x,r))}.
\end{equation*}
\end{enumerate}
\end{example}

\begin{remark}
The above example includes the case where $X$ is a Lebesgue space, a
Lorentz space, a Lorentz-Zygmund space or a generalized
Lorentz-Zygmund spaces.
\end{remark}

\begin{remark}
In the particular case $\Omega =\mathbb{R}^{n},$ $E=L^{q}$ ($1\leq
q\leq \infty )$, $b=1$ and $X=L^{p}$ $(1\leq p<\infty )$ we obtain
the following well-known result:

\begin{enumerate}
\item  If $\frac{1}{p}>\frac{\theta }{n},$ then
\begin{equation*}
\left( \int_{0}^{\infty }\left( f^{\ast \ast }\left( t\right) t^{\frac{1}{p}-%
\frac{\theta }{n}}\right) ^{q}\frac{dt}{t}\right) ^{1/q}\preceq
\left\| f\right\| _{B_{p,q}^{\theta }}.
\end{equation*}

\item  If $\frac{\theta }{n}=\frac{1}{p},$ then
\begin{equation*}
\left\| f\right\| _{L^{\infty }}\preceq \left\| f\right\|
_{B_{p,1}^{\theta }},
\end{equation*}
and $f$ is and essentially continuous.$\ $
\begin{equation*}
\left( \int_{0}^{1}\left( \frac{f^{\ast \ast }\left( s\right) }{1+\ln \frac{1%
}{s}}\right) ^{q}\frac{ds}{s}\right) ^{1/q}\preceq \left\| f\right\|
_{B_{p,q}^{\theta }}+\left\| f\right\| _{L^{1}+L^{\infty }}\text{ \ (if }%
1<q<\infty \text{),}
\end{equation*}
and
\begin{equation*}
\left\| f\right\| _{bmo}\preceq \left\| f\right\| _{_{B_{p,\infty
}^{\theta }}}.
\end{equation*}

\item  If $\frac{1}{p}<\frac{\theta }{n},$ then
\begin{equation*}
\left\| f\right\| _{L^{\infty }}\preceq \left\| f\right\|
_{B_{p,q}^{\theta }}
\end{equation*}
and
\begin{equation*}
\left| f(x)-f(y)\right| \preceq \left\| f\right\| _{B_{p,q}^{\theta
}(\Omega )}d(x,y)^{\theta -\frac{n}{p}}.
\end{equation*}
\end{enumerate}
\end{remark}

\section{Appendix}

In this section we give the proof of the lemma \ref{Suerte}, which
follows from a persual of the proofs of Lemma 2.3 and Proposition
2.5 of \cite{Ran1} combined with H\"{o}lder's inequality
(\ref{Hold}).

\begin{proof}
Part 1) By H\"{o}lder's inequality, we have
\begin{align*}
\left( T_{R}^{1}f(x)\right) ^{q}& :=\left( \dint_{B(x,R)}\left|
f(y)\right|
d\mu (y)\right) ^{q} \\
& \leq \left( \dint_{B(x,R)}\left| f(y)\right| ^{p}d\mu (y)\right) ^{q/p} \\
& =\mu (B(x,R))^{q/p}\int_{B(x,R)}\left| f(y)\right| ^{p}d\mu
(y)\left(
\int_{B(x,R)}\left| f(y)\right| ^{p}d\mu (y)\right) ^{-1+q/p} \\
& \leq \mu (B(x,R))^{q/p}\int_{B(x,R)}\left| f(y)\right| ^{p}d\mu
(y)\left( \left\| \left| f\right| ^{p}\right\| _{X}^{-1+q/p}\varphi
_{X^{\prime
}}\left( \mu (B(x,R))\right) ^{-1+q/p}\right) \\
& \leq \frac{\left\| \left| f\right| ^{p}\right\|
_{X}^{-1+q/p}}{\inf_{x\in \Omega }\varphi
_{X}(B(x,R))^{q/p-1}}\left( \dint_{B(x,R)}\left| f(y)\right|
^{p}d\mu (y)\right) \text{ (by (\ref{si}))}.
\end{align*}
Taking $\left\| \cdot \right\| _{X}$%
\begin{align*}
\left\| \left( T_{R}^{1}f(x)\right) ^{q}\right\| _{X}& \leq
\frac{\left\| \left| f\right| ^{p}\right\| _{X}^{-1+q/p}}{\inf_{x\in
\Omega }\varphi _{X}(B(x,R))^{q/p-1}}\left\| \dint_{B(x,R)}\left|
f(y)\right| ^{p}d\mu
(y)\right\| _{X} \\
& \preceq \frac{\left\| \left| f\right| ^{p}\right\| _{X}^{-1+q/p}}{%
\inf_{x\in \Omega }\varphi _{X}(B(x,R))^{q/p-1}}\left\| \left|
f\right|
^{p}\right\| _{X}\text{ (by Lemma \ref{Suerte} (i))} \\
& =\frac{\left\| \left| f\right| ^{p}\right\| _{X}^{q/p}}{\inf_{x\in
\Omega
}\varphi _{X}(B(x,R))^{q/p-1}} \\
& \preceq \frac{\left\| f\right\| _{X^{(p)}}^{q}}{\varphi _{X}(\min
(R^{k},R^{n}))^{q/p-1}}\ \text{(by (\ref{belw})).}
\end{align*}

Part 2) Let us write
\begin{equation*}
I_{R}(f,x)^{q}=\left( \dint_{B(x,2R)}\dint_{B(x,R)}\left|
f(y)-f(z)\right| d\mu (z)d\mu (y)\right) ^{q}
\end{equation*}
and
\begin{equation*}
\nabla _{R}^{p}f(x)=\left( \dint_{B(x,R)}\left| f(y)-f(x)\right|
^{p}d\mu (y)\right) ^{1/p}.
\end{equation*}
From H\"{o}lder's inequality we can deduce that (see \cite[Page
7]{Ran1})
\begin{equation*}
I_{R}(f,x)^{q}\preceq \frac{1}{\mu (B(x,R))^{q/p-1}}\left(
\int_{B(x,3R)}\dint_{B(x,3R)}\left| f(y)-f(z)\right| ^{p}d\mu
(z)d\mu (y)\right) ^{-1+q/p}\left( \nabla _{2R}^{p}f(x)\right) ^{p}.
\end{equation*}
On the other hand by (\ref{Hold}),
\begin{align*}
& \left( \int_{B(x,3R)}\dint_{B(x,3R)}\left| f(y)-f(z)\right|
^{p}d\mu
(z)d\mu (y)\right) ^{-1+q/p} \\
& \leq \left\| \dint_{B(x,3R)}\left| f(y)-f(z)\right| ^{p}d\mu
(z)\right\|
_{X}^{-1+q/p}\varphi _{X^{\prime }}\left( \mu (B(x,3R))\right) ^{-1+q/p} \\
& =\left\| \left( \nabla _{3R}^{p}f(x)\right) ^{p}\right\|
_{X}^{-1+q/p}\varphi _{X^{\prime }}\left( \mu (B(x,3R))\right)
^{-1+q/p}.
\end{align*}
So, using the concavity of $\varphi _{X}$ and the fact that $\mu $
is a doubling, we get
\begin{align*}
I_{R}(f,x)^{q}& \preceq \left( \frac{\varphi _{X^{\prime }}\left(
\mu (B(x,3R))\right) }{\mu (B(x,R))}\right) ^{q/p-1}\left\| \left(
\nabla _{3R}^{p}f(x)\right) ^{p}\right\| _{X}^{-1+q/p}\left( \nabla
_{2R}^{p}f(x)\right) ^{p} \\
& \preceq \frac{1}{\inf_{x\in \Omega }\varphi
_{X}(B(x,R))^{q/p-1}}\left\| \left( \nabla _{3R}^{p}f(x)\right)
^{p}\right\| _{X}^{-1+q/p}\left( \nabla _{3R}^{p}f(x)\right) ^{p},
\end{align*}
and, taking the $X$ norm in the above expression we get
\begin{equation*}
\left\| I_{R}(f,x)^{q}\right\| _{X}\preceq \frac{1}{\inf_{x\in
\Omega }\varphi _{X}(B(x,R))^{q/p-1}}\left\| \left( \nabla
_{3R}^{p}f(x)\right) ^{p}\right\| _{X}^{q/p},
\end{equation*}
or equivalently
\begin{equation}
\left\| I_{R}(f,x)\right\| _{X^{(q)}}\preceq \frac{1}{\inf_{x\in
\Omega }\varphi _{X}(B(x,R))^{1/p-1/q}}\mathcal{E}_{X^{(p)}}(f,3R).
\label{suma0}
\end{equation}
On the other hand, for every integer $m,$ we have
\begin{align}
& \dint_{B(x,R)}\dint_{B(x,R/2^{m+1})}\left| f(y)-f(z)\right| d\mu
(z)d\mu
(y)  \label{suma} \\
& \leq
\sum_{i=0}^{m}\dint_{B(x,R/2^{i})}\dint_{B(x,R/2^{i+1}))}\left|
f(y)-f(z)\right| d\mu (z)d\mu (y).  \notag
\end{align}
Consequently,
\begin{align*}
& \left\| \dint_{B(x,R)}\dint_{B(x,R/2^{m+1})}\left|
f(y)-f(z)\right| d\mu
(z)d\mu (y)\right\| _{X^{(q)}} \\
& \leq \sum_{i=0}^{m}\left\|
\dint_{B(x,R/2^{i})}\dint_{B(x,R/2^{i+1}))}\left| f(y)-f(z)\right|
d\mu (z)d\mu (y)\right\| _{X^{(q)}}.
\end{align*}
On the other hand, for every non-negative integer $m$ we have
\begin{align}
& \left\| \dint_{B(x,R)}\left| f(y)-\dint_{B(x,R/2^{m+1})}f(z)d\mu
(z)\right| d\mu (y)\right\| _{X^{(q)}}  \label{zzz} \\
& \leq \left\| \dint_{B(x,R)}\dint_{B(x,R/2^{m+1})}\left|
f(y)-f(z)\right| d\mu (z)d\mu (y)\right\| _{X^{(q)}}.  \notag
\end{align}
From Lebesgue's differentiation theorem (for the doubling measures),
Fatou's lemma and letting $m\rightarrow \infty $, it follows that
\begin{align*}
E_{X^{(q)}}(f,R)& =\left\| \dint_{B(x,R)}\left| f(y)-f(x)\right|
d\mu (z)d\mu
(y)\right\| _{X^{(q)}} \\
& =\left\| \dint_{B(x,R)}\lim_{m\rightarrow \infty }\left|
f(y)-\dint_{B(x,R/2^{m+1})}f(z)d\mu (z)\right| d\mu (y)\right\| _{X^{(q)}} \\
& \leq \left\| \lim_{m\rightarrow \infty
}\dint_{B(x,R)}\dint_{B(x,R/2^{m+1})}\left| f(y)-f(z)\right| d\mu
(z)d\mu
(y)\right\| _{X^{(q)}}\text{ \ (by (\ref{zzz}))} \\
& \leq \left\| \sum_{i=0}^{\infty
}\dint_{B(x,R/2^{i})}\dint_{B(x,R/2^{i+1}))}\left| f(y)-f(z)\right|
d\mu
(z)d\mu (y)\right\| _{X^{(q)}}\text{ \ (by (\ref{suma}))} \\
& \leq \sum_{i=0}^{\infty }\frac{1}{\inf_{x\in \Omega }\varphi
_{X}(B(x,R/2^{i+1}))^{1/p-1/q}}\mathcal{E}_{X^{(p)}}(f,3R/2^{i+1})\text{
\
(by (\ref{suma0}))} \\
& \leq \sum_{i=0}^{\infty }\frac{1}{\inf_{x\in \Omega }\varphi
_{X}(B(x,R/2^{i-2}))^{1/p-1/q}}\mathcal{E}_{X^{(p)}}(f,3R/2^{i-1}) \\
& \leq \sum_{i=0}^{\infty }\int_{R/2^{i-1}}^{R/2^{i-2}}\frac{\mathcal{E}%
_{X^{(p)}}(f,r)}{\inf_{x\in \Omega }\varphi _{X}(B(x,r))^{1/p-1/q}}\frac{dr}{%
r}\text{ \ (by (\ref{qc}))} \\
& \preceq \int_{0}^{4R}\frac{\mathcal{E}_{X^{(p)}}(f,r)}{\varphi
_{X}\left( \min (r^{k},r^{n})\right) ^{1/p-1/q}}\frac{dr}{r}\text{ \
(by (\ref{belw})).}
\end{align*}
\end{proof}

\section*{Conflict of interest} The authors declare that they have no conflict of interest.

\section*{Data availability statement}
Data sharing not applicable to this article as no datasets were
generated or analyzed during the current study.

\end{document}